\pdfoutput=1
\documentclass[11pt]{article}
\usepackage[a4paper,margin=2.5cm]{geometry}
\usepackage[english]{babel}
\usepackage{amsmath}
\usepackage{amssymb}
\usepackage{graphicx}
\usepackage{hyperref}
\usepackage[capitalise,nameinlink,noabbrev]{cleveref}

\usepackage{setspace}
\usepackage{xspace}
\usepackage{float}

\usepackage{booktabs}

\usepackage[numbers]{natbib}

\title{New Heuristics for Large-Scale Joint Order Selection,\\ Allocation, Batching, and Picking}
\author{%
  Alexander Grosz, Leander Schnaars\\[0.5em]
  \small Chair of Operations Research\\[-0.2em]
  \small Technical University of Munich\\[-0.2em]
  \small Arcisstraße 21, 80333 Munich, Germany\\[-0.2em]
  \small \texttt{alexander.grosz@tum.de}, \texttt{leander.schnaars@tum.de}
}
\date{}

\newcommand*{\N}{\mathbb{N}} \newcommand*{\Z}{\mathbb{Z}} \newcommand*{\Articles}{\mathcal{A}}
\newcommand*{\Items}{\mathcal{I}} \newcommand*{\Orders}{\mathcal{O}} \newcommand*{\Zones}{\mathcal{Z}}
\newcommand*{\Capacity}{\mathcal{C}} \newcommand*{\Goal}{\mathcal{G}} \newcommand*{\Batch}{\mathcal{B}}
\newcommand*{\Feasible}{\mathcal{F}}
\newcommand*{\card}[1]{\left|#1\right|}

\newcommand*{\algo}[1]{\texttt{#1}\xspace}

\newcommand*{\dga}{\algo{DGA}}
\newcommand*{\fdga}{\algo{FDGA}}
\newcommand*{\fdgaskip}{\algo{FDGA-Skip}}
\newcommand*{\rdga}{\algo{RDGA}}
\newcommand*{\rdgan}[1]{\algo{RDGA-#1}}
\newcommand*{\rdgak}{\rdgan{$k$}}
\newcommand*{\rdgasqrt}{\rdgan{$\sqrt{|\Orders|}$}}
\newcommand*{\ldga}{\algo{LocalDGA}}

\newcommand*{\fdgasoft}{\algo{FDGA-Soft}}
\newcommand*{\fdgasoftg}[1]{\algo{FDGA-Soft(#1)}}
\newcommand*{\fdgalinear}{\algo{FDGA-Linear}}
\newcommand*{\fdgalinearn}[1]{\algo{FDGA-Linear(#1)}}
\newcommand*{\fdgacutoff}{\algo{FDGA-Cutoff}}
\newcommand*{\fdgacutoffn}[1]{\algo{FDGA-Cutoff(#1)}}

\begin{document}

\newcommand{\TABLE}[3]{%
  \refstepcounter{table}%
  \begin{center}
    \small
    \textbf{Table \thetable.} #1\par
    \vspace{0.75em}
    #2\par
    \vspace{0.75em}
    \parbox{\linewidth}{#3}%
  \end{center}%
}
\newcommand{\FIGURE}[3]{%
  \refstepcounter{figure}%
  \begin{center}
    #1\par
    \textbf{Figure \thefigure.} #2\par
    \smallskip
    \parbox{\linewidth}{\footnotesize #3}%
  \end{center}%
}
\newenvironment{APPENDICES}{\appendix}{}

\maketitle

\begin{abstract}
  Manual order picking in warehouses plays an important role in modern e-commerce fulfillment.
  In this work, we consider a warehouse fulfillment model recently introduced by researchers at Zalando, a leading fashion e-commerce company.
  It simultaneously models the hierarchy of decision stages in a very large-scale setting with a semi-online order set, where only a certain number of items have to be fulfilled in each time interval.
  For their base model, we present both new and improved algorithms together with comprehensive performance evaluation.
  We further extend the model to a setting with multiple phases, which more closely represents the situation observed in real warehouses, and adapt the algorithms to this setting, including new variants that incorporate multi-phase specific features.
  We demonstrate that this extension allows us to model, control, and optimize for certain parameters such as order age which play an important role in managerial considerations with the growing importance of next or even same day delivery.
\end{abstract}

\noindent\textbf{Keywords:} Large-scale optimization, warehouse fulfillment, order selection and picking

\section{Introduction}\label{sec:intro}
In the age of ever-growing e-commerce, manual order picking in large-scale warehouses still plays an important part in the fulfillment supply chain. While semi- or fully-automated warehouses are seeing increasing use, manual order picking is still seen in a majority of warehouses \citep{Boysen2019-WarehousingInECommerce}. In a typical setting, items are distributed on shelves throughout the warehouse and pickers collect them individually. The cost of order picking is often cited to make up more than 50\% of the total cost of the warehouse operation. Pickers, on average, spend around half of their time traveling between their assigned items to pick or the drop-off stations \citep{Tompkins2010-FacilityPlanning}. This implies that there are major cost-saving opportunities just through optimization of this part of the supply chain. In this work we focus on a setting recently introduced by researchers working at Zalando, one of the largest e-commerce fashion retailers in Europe \citep{ZalandoPaper}.
Their setting combines a multi-layered decision process with a very large-scale setting, with tens of thousands of pending orders at any given point in time.

\subsection{Warehouse Order Picking}
A multitude of variants of order picking problems in real warehouses exist and an even larger amount of settings has been studied from a theoretical perspective. Several recent surveys provide a comprehensive overview over the current state of the art \citep{Boysen2019-WarehousingInECommerce,Masae2020-OrderPickingLitRev,Pinto2022-BatchingReview,CVMTB:TrendsInOrderPicking,Bock2024-WarehouseSurvey}. While some publications discuss the optimization of the warehouse layout itself, most works consider a fixed layout which has been determined during the warehouse's initial design phase. Commonly studied layouts are rectangular shelf arrangements with one or more cross-aisles, fishbone configurations and V-shaped rack configurations \citep{Roodbergen2006-warehousedesign,Cardona2015-warehousedesign,Zhou2019-warehousedesign,Bocelli2024-warehousedesign}. Items are spread across the shelves throughout the warehouse. Unlike in settings such as spare parts management for big industrial machines, in e-commerce often many items of a specific article (SKU) are stocked inside of the warehouse. It has been shown to be beneficial to spread these items throughout the warehouse instead of keeping all at the same location, as this enables routing synergies when collecting items \citep{Onal2017-explosive}, so modern warehouses typically support such storage policies. This strategy is called scattered or explosive storage.

In our setting, we focus on a pure picker-to-parts model, where humans, possibly with the help of some form of cart, walk throughout the warehouse and collect the items. Some other variants have been studied where the pickers are assisted by technological tools such as automated guided vehicles which deliver picked orders to the depots or rideable ``cobots'' which enable pickers to travel across large distances more quickly \citep{Loffler2022-AGVPicking,Tutam2024-cobot}. The pending demand is usually specified by a list of orders, each of which consists of one or several articles contained in the warehouse. While some low throughput systems exist where each order is individually collected by a single picker, in larger warehouses with high order volume the task is typically spread amongst several pickers who each pick articles belonging to different orders. For each picker, an efficient route from the depots through the warehouse has to be calculated. The picked items are then combined into the respective orders in a later step of the logistics process, which we do not consider here. This setting and closely similar ones have been studied by many authors. Focusing on the picklist route step, results in \cite{Ratliff1983,Roodbergen2001,Cambazard2018-np,Pansart2018} have shown polynomial time solvability of this subproblem for certain restricted warehouse layouts. Settings which combine batching/assignment and routing have been studied in works such as \cite{Briant2020-routeassignment,Aerts2021-routeassignment,Tran2024-routeassignment} and are usually computationally hard.

A novelty in the model of \cite{ZalandoPaper}, which we base our work on, is the introduction of an \emph{item goal}, such that not all orders have to be fulfilled, but rather a subset of orders containing a set number of items has to be picked. This enables additional synergies to be exploited, as not all orders have to be picked in the considered step, so orders with complementing article placement can be preferred, while the remaining orders can benefit from new orders arriving later.

\subsection{Large-Scale Warehouses}
The model considered in \cite{ZalandoPaper} and here focuses on a very large-scale setting as faced by the largest players in modern e-commerce. Zalando itself has more than $1.8$ million different products available \citep{ZalandoReport2024} and thus the considered scales vastly exceed most of the order picking literature. While not every one of Zalando's warehouses contains all of the products, we model several hundred thousand distinct articles with a million total items contained in our instances spread across $100$ zones with $100 \times 100$ locations each. For our largest instances, around 50,000 orders are considered at any time, of which 10,000 need to be batched. Most previous works consider a few hundred orders at most, with single zone warehouses and a small number of pickers and articles.
A review of recent relevant articles shows typical order counts of a few hundred up to a thousand, with at most a few thousand item positions in the warehouse. In \cite{Rahmani_Mokarrari2025-ordersize}, the authors consider a setting with around 10,000 item locations and at most $700$ orders. For leaf layout warehouses, \cite{Masae2021-ordersizeleaf} review common picker routing algorithms and evaluate their performance for a single zone and up to 1,000 orders. Further numerical evaluations with between $600$ and 7,200 item locations and up to $100$ and $200$ orders can be found in \cite{Kubler2020-ordersize1,Cao2023-ordersize2,Cano2021-ordersize3,Aerts2021-routeassignment,Aksu2025-ordersize5}.

\subsection{Semi-Online Time Scales}
The aim of the item goal in this model is to represent a real-time semi-online setting, where in some discrete steps of time, usually between 5 and 60 minutes, the current state of the pending orders is considered and then an order batch and subsequent picklist assignment is created. Compared to a fully online setting, this establishes a precomputed schedule for the respective time period, which can increase reliability. Due to this, the overall algorithm has to run fairly quickly, with target times of a few minutes at most. Due to the increased offering of services such as same day shipping and express delivery \citep{Fisher2019-deliveryadvantage,Cui2024-deliveryadvantage}, it becomes increasingly important to dynamically update the order set considered by the algorithm timely, as otherwise the targeted quick delivery times could not be achieved and possible synergies when picking the orders would be missed. Most relevant literature in order picking views the problem in a more classical static setting with a fixed set of orders and a fixed warehouse inventory, searching for an optimal or close to optimal solution to the respective problem configuration. Such approaches generally do not achieve the restrictive target times for such a dynamic setting and faster but more inaccurate heuristics have to be employed.

\subsection{Hierarchy of Subproblems}

As mentioned above, there exists a plethora of specialized problems in warehouse optimization and all of them interact in some way or another with whichever optimization objective is being chosen. Thus, many (especially more recent) publications in the field focus on integrating several of these problems into a single optimization problem to achieve the best outcome. Studying these problems in isolated settings often reveals that they are already hard from the perspective of complexity theory and their combination reaches limits of feasibility in the real-world applicability of the solutions very quickly, especially considering the size of real-world instances and the necessity of providing solutions reasonably quickly. For extensive discussions on both the individual levels, as well as various integrations of them, we refer to recent surveys on this field \citep{PGAD:OrderBatchingReview, Masae2020-OrderPickingLitRev, CVMTB:TrendsInOrderPicking}.
In the following, we present the relevant hierarchical stages of the problem we discuss in this work. These are the stages we attempt to solve simultaneously as a single integrated problem.

\paragraph{Order selection}

Typically, the full set of orders available to the problem is expected to be completed at once, and the optimization procedure attempts to find an optimal solution to do so. In the model considered in this paper, we allow this requirement to be relaxed to some degree: from the current set of orders, we only attempt to complete a certain number (given by what we call the \emph{item goal}) of them, while the remainder is postponed to a later time period. The idea behind this new approach compared to the existing literature is that new orders coming in during a later stage could provide opportunity for better synergies with the set of remaining orders.

The time period modeled in this step represents between tens of minutes to a few hours of real-world time, after which a significant amount of new orders have already been placed within the scale of the application we model here. This additionally addresses the situation that customers are able to place orders at any time, while warehouses often have limited working hours. The model can address these imbalances in the size of the full set of orders while maintaining consistent warehouse throughput.

\paragraph{Order batching}

The discretization we choose in the previous step works on a larger scale than what is required for the operational scale in the warehouse. On smaller time intervals, we need to provide the subsequent stages of the order fulfillment process with all articles required to complete a certain set of orders. Clearly, some maximum throughput capacity sets the limit to how many orders should be fetched simultaneously. We control this selection process by grouping orders into batches. All respective items of the selected orders in a single batch are then retrieved from the warehouse in parallel by multiple pickers across multiple warehouse zones and then forwarded to the next stage before continuing with the next batch of orders.
This approach not only allows a single picker to be assigned items for multiple orders to be retrieved within a single tour (according to the capacity of their cart), but also a single order can be spread across different pickers to improve the efficiency of the process.

\paragraph{Item allocation}

With the introduction of mixed-shelves storage and the explosive storage policy, where identical items of the same article do not have to be stored at a single warehouse location, another decision has to be considered in the optimization process: which specific item will be picked to fulfill the article request of a given order. While this additional step introduces significant amounts of computational complexity into the problem, the synergistic opportunity of this decision during the routing step is similarly significant.

\paragraph{Picker routing}

Finally, with the set of orders selected, partitioned into batches, and their respective articles being assigned a specific item in the warehouse, the pickers are sent out to retrieve these items from the warehouse. We consider a model where the warehouse is split into separate zones to which the pickers are restricted for a single picklist. Within a batch, the solution consists of a set of picklists, describing for every picker in which order the items should be fetched from the warehouse such that the total travel distance is minimized. Again, we emphasize that while this last optimization step is in itself a mathematically intricate optimization problem, its optimal solution is highly affected by all the choices higher up in the hierarchy presented here, and therefore an integrated approach is favorable.

\subsection{Contributions}
In their article, \cite{ZalandoPaper} introduced a new setting for order picking in warehouses, inspired by their knowledge about the warehouse operations of Zalando. In this work, we build upon their model and algorithms and thus provide the next step in exploring both algorithms for and limitations of the model. To make the model more representative of real world constraints and the analysis less susceptible to local optima, we extend it to a multi-phase setting and analyze all algorithms in this new setting. The analysis shows that our algorithms both run in reasonable time while also providing good heuristic cost results. Additionally, some variants also achieve strong performance with regards to some secondary managerial objectives like fulfillment guarantees and delay prevention. The main contributions of this paper are summarized as follows:
\begin{itemize}
  \item We present an improved and computationally efficient implementation of a reference algorithm for the joint order selection, allocation, batching, and picking problem.
  \item Several new algorithmic variants for this problem with attractive runtime behavior maintain good solution quality in computational experiments.
  \item We conduct empirical analysis of these algorithms on synthetic data and discuss structural properties of the solutions.
  \item An extension of the model to a multi-phase setting allows for more accurate analysis of the algorithms' performances.
  \item We introduce several novel age-dependent algorithms to address shortcomings of the established algorithms over multiple phases, especially with regard to order starvation and fulfillment guarantees which are important from a managerial perspective.
  \item Through a broad empirical analysis of the algorithms on the multi-phase model, we confirm the desired long-term behavior.
\end{itemize}

Our implementation of the models and algorithms can be found at \url{https://github.com/puenka/multiphase-batching-heuristics}.

\section{Problem Definition}

Our paper is based on the problem and model presented by researchers from Zalando in \cite{ZalandoPaper}, which we
paraphrase in this section and extend later. Let us first discuss the underlying logistics of the problem on a less mathematical level
from the bottom to the top of the hierarchy.
An instance of this problem will be called a \emph{phase} in the multi-phase model introduced later in \cref{sec:multiphase-model}.

\begin{figure}
  \centering
  \includegraphics[scale=0.7]{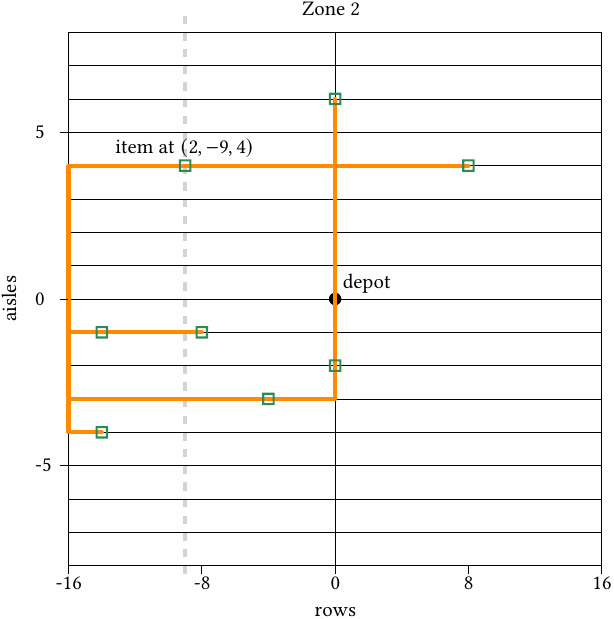}
  \caption{An exemplary warehouse setup, indicating a route that picks up several items and an item with its location indices (the corresponding row is highlighted by a dashed line).}
  \label{fig:warehouselayout}
\end{figure}

The warehouse considered in the problem is split into independent zones, in which the \emph{picking} part of the problem
is located. Each zone is a physical warehouse section with locations indicated by racks and aisles. At each of these
locations, individual items are located and these locations are therefore identified by the zone, rack, and aisle
number. The picker can navigate along aisles and change between different aisles at three orthogonal cross-aisles at the
center and at both ends of each aisle. Each route for the picker starts at a depot located in the
center of the zone and they have to return to this depot at the end of the route. Also, the routes cannot span across
different zones, so after the item selection the routing problem is limited to a single zone. For simplicity, we
typically assume a symmetric warehouse setup. See \cref{fig:warehouselayout} for an illustration of a warehouse zone with several items being picked up in a route.

The next level in the hierarchy is the \emph{item selection}: For each article (note the differentiation between
\emph{article} and \emph{item}) that has to be picked up by the picker to complete an order, one of possibly several
identical items at different locations in the warehouse has to be chosen.
Additionally, each article has an associated \emph{volume}, which reflects the portion it takes up in the bounded
capacity of a single picker as a single numerical value.

The main objective is for a given list of \emph{orders} to be completed. Each order consists of a list of articles and
orders are not allowed to be fulfilled only partially. The \emph{order selection} part of the problem consists of
selecting a given number of orders in each \emph{batch}. The maximum number of orders in each batch reflects the
capacity of the warehouse and the processing capacity after the item collection in the logistics chain for a given time
frame.

Furthermore, a novelty in \cite{ZalandoPaper} was the introduction of an \emph{item goal}. Instead of requiring all
orders in the problem to be completed, like it is observed in most other picker routing problems, we focus on simply achieving a set number of items in each phase of the problem, while keeping a backlog of orders to be fulfilled in a later phase.
This potentially enables additional optimization potential, as we can select a subset of orders to be completed that
synergize well with each other and hope for new orders to come in during later phases, that now complement the remaining orders better.

\subsection{Mathematical Model} \label{sec:singlephase-model}

We briefly review the most important parameters of the model here, as some of them are used in construction and specification of the algorithms, for full details regarding all parameters and the associated constraints and goals see Appendix \ref{sec:app:mathmodel}. For ease of notation we use the term \emph{set} to refer to multisets.

During execution of any of the algorithms in this work, the state of the warehouse is captured by several variables. Denote by $\Articles$ the set of articles and by $\Items$ the set of items in the warehouse. Recall that articles correspond to stock keeping units, while items correspond to the actual physical objects located in the warehouse. The volume of article $a\in \Articles$ is given by $v_a \in \N$. Each item $i\in \Items$ can fulfill its corresponding article $a_i \in \Articles$ and is located at a location triple $(z_i, r_i, s_i)\in \Z^3$ corresponding to zone, rack, and aisle number. During execution of the algorithms, the set of unfulfilled orders is given by $\Orders$, where each order consists of a set of articles. Each of the pickers in the model has a capacity $\Capacity \in \N$, which captures how much volume they can transport at once.

\subsection{Model Extension to Multiple Phases}
\label{sec:multiphase-model}
The multi-phase model is a direct extension of the single-phase model to a larger timeframe. It can be seen as a sequence of single phases, with restocking periods in-between. This is meant to model the continuous process of restocking and order fulfillment periods more closely, instead of only considering a single isolated phase.

After each phase, the warehouse is restocked with articles and the list of remaining orders is  refilled. Both steps are done in a way to ensure that the number of items in the warehouse remains constant throughout the whole model run. Additionally, through appropriate design of the newly added orders, it is ensured that every order is theoretically fulfillable at any given point in time, so there are never orders which contain articles not contained in the warehouse.

Note that iterative application of any algorithm for the single-phase model gives solutions for multiple phases as well. We investigate the structure of these solutions and address their issues in \cref{sec:multiphase-empirical-analysis}.

\section{Complexity and Polynomial Solvability}
In this section, we briefly review the state of the art regarding complexity results and solvability of picking problems. We additionally discuss further sub-problems related to the multi-zone setting and to picker capacities. Such analysis is important to understand which approaches and algorithms are reasonable for the given settings of order picking.

\subsection{General Order Picking Complexity}
The complexity of general order picking problems is well studied, with most problems shown to be NP-hard, unless heavily restricted.
This is due to order picking being interpretable as a combination of several classical well-known hard problems, such as Traveling Salesperson (TSP), Steiner-tree, and Bin-packing \citep{Karp1972-21probs}.

The setting of so called wave-picking, where orders are released in large batches and the next batch can only be started when the previous one is completed, was studied in \cite{Gademann2001-cv}, with NP-hardness established there. For an order batching subproblem used to reduce picker traveling times, \cite{Gademann2005-jg} show hardness if some batches contain more than two items. In their work, \cite{Schiffer2022-ij} show that a version they call Generalized Order Picking Problem is NP-hard, even if all picklists are restricted to be of length one and the warehouse has a simple structure with a single zone and only two cross-aisles. Focusing mainly on the TSP aspect of order picking, in their survey \cite{Bock2024-WarehouseSurvey} present a wide variety of hardness proofs and some polynomial algorithms for many variants of the problem.

\subsection{Efficiently Solvable Instances}
\label{sec:efficientrouting}
In contrast to the NP-hardness of general order picking problems, many restricted variants or subproblems are known to be solvable in polynomial time. By assuming certain parameters such as the number of aisles and rows to be of a fixed size, parametrized polynomial time algorithms are also known.

Order picking in a rectangular warehouse with only two cross-aisles at the two ends of the aisles was shown to be solvable in linear time in the number of aisles by \cite{Ratliff1983}. Recently, the runtime for these two cases was improved in \cite{dunn2025det}. This setting was extended by \cite{Roodbergen2001} with the introduction of additional cross-aisles in the inner part of the warehouse at the cost of a much higher runtime. For an arbitrary fixed number $h$ of cross-aisles, \cite{Cambazard2018-np} introduced a polynomial time algorithm with an exponential dependency on $h$. Additional analysis and further details can be found in \cite{Pansart2018}. A more general framework which encompasses these types of settings together with some optimizations which slightly improve the runtime can be found in \cite{Schiffer2022-ij}.

\subsection{NP-Hardness of Order Selection}
The complexity and solvability of the picker routing problem is well studied, as demonstrated in the previous two sections. In our setting, with multiple disjoint independent zones and articles of which copies can be spread throughout these zones, additional complexities arise. As each order consists of several articles, each of which might have several corresponding warehouse items, we have to decide which of them to pick to fulfill the order.

For a given set of orders, assigning articles to orders to fulfill all demands without assigning an article twice is a bipartite matching problem, which is well known to be solvable in polynomial time. However, if we also want to optimize certain cost parameters at the same time, the problem becomes much harder. We call a selection of articles \emph{satisfying} for a set of orders if all item demands from these orders are fulfilled by the set. We claim that deciding whether a satisfying set of articles which can be picked for a certain cost exists is NP-hard to decide, even when the warehouse structure is very simple and optimal picking routes are known. The same hardness extends to deciding whether a satisfying set exists which only uses a certain number of zones. There are multiple ways to show this, which can be used to establish the hardness even for more restricted instances.

Even in the case where there is only one unique article corresponding to each item in the warehouse, certain subproblems remain difficult. In this setting, for each set of orders there is a unique set of satisfying articles. Even if these items were all situated in the same zone at the same place, due to their volumes and the maximum capacity of each picker, it would be hard to decide whether there is a way to pick the items using two pickers or whether at least three are required. This follows as a direct consequence of the hardness of the classic Bin-packing problem.

\section{Single Phase Setting}\label{sec:singlephase}
In this section we introduce and evaluate several algorithms for the single-phase setting. We focus on high level descriptions and refer to the implementation for finer details. The Distance Greedy Algorithm (\dga) and the base version of the Randomized Distance Greedy Algorithm (\rdga) serve as reference algorithms and are only slightly adapted from \cite{ZalandoPaper}. We present additional variants and improved versions of these algorithms. The remaining algorithms are novel and use a variety of approaches. We first introduce all algorithms and then give a complete analysis and comparison of their performance in \cref{sec:singlephase-empirical-analysis}. In \cref{sec:singlephase-performance} we examine and compare the performance impact of certain subroutine choices, such as the optimization of the picklist routes by better solvers.

\subsection{Distance Greedy Algorithm}
The \emph{Distance Greedy Algorithm} (\dga) is logically identical to the one in \cite{ZalandoPaper}. For this algorithm, the item batches are created sequentially and during the creation of each batch, one order after another is picked, until the maximum size of the batch is reached. The respective next order is picked based on a distance function which estimates the cost to pick each order based on the previously selected orders and thus prioritizes orders whose items are close to already selected items or close to the depots. The picklist creation is done at the end through a simple greedy procedure. In \cref{sec:singlephase-performance}, we empirically compare this greedy procedure to an improved solution computed by modeling this subproblem as a constrained vehicle routing problem (CVRP) using a solver from the external library OR-Tools \citep{ortools_routing}.

As noted by the authors in \cite{ZalandoPaper} and provided in their codebase, in a naive implementation, this algorithm runs too slowly on larger instances, with runtimes reaching several hours. Through various design and implementation improvements (all within a pure Python implementation to maintain comparability to the reference), we have achieved a version of this algorithm with equivalent computational output, but vastly improved runtime performance, being fast enough to be employed in a setting where results within minutes are required. The most important change is keeping track of already computed distance estimates within a heap data structure and updating them only when necessary. After selecting the next order and its associated items, using several structural insights into efficient picker tours, the algorithm analyzes which items and thus orders could potentially have their cost estimates affected by the order choice, which leads to much fewer required updates and a resulting significantly faster runtime. For easier disambiguation and comparison we call our improved version \emph{fast} \dga (\fdga).

We also introduce a variant of this algorithm based on a simplification of the subroutine which selects the next order. To estimate the marginal cost of an order, each article is iteratively assigned the currently best fitting item. While \dga and \fdga then compare the addition of subsequent items also to the already picked items for the order under investigation, the variant \fdgaskip (which we adapt from the improved implementation \fdga) omits this part of the computation and only estimates the cost by comparison to all items of previously selected orders.

\subsubsection{Randomized Distance Greedy Algorithm}
The base version of the \emph{Randomized Distance Greedy Algorithm} (\rdga) from \cite{ZalandoPaper} is a severe simplification of the \dga. Instead of using the cost estimate to select the next order, a random unfulfilled order is picked. Then, items are selected that fulfill this order which heuristically minimize the marginal increase of the picklist cost with regard to the already selected orders.

We build upon this idea of randomization by introducing a sampling-based variant called \rdgak. Our goal is to achieve a trade-off between the improved runtime obtained by considering only a subset of the remaining orders and the increased resulting total cost through lost synergy potential.
Instead of picking a single order like \rdga, our new variant \rdgak uniformly samples $k$ of the remaining orders and then selects the one minimizing the expected cost. We also allow $k$ to be a function depending on the state of the warehouse, such as the number of remaining orders, or more specifically, their square root, called \rdgasqrt.

\subsubsection{Local Distance Greedy Algorithm}
Our new algorithm called Local Distance Greedy Algorithm (\ldga) attempts to refine the sampling approach of the \rdga. In this algorithm, we aim to preserve the improved runtime obtained by omitting the evaluation of a large number of potential orders, while reducing the loss in optimality incurred through the random choices.
Instead of resampling the entire space of orders as is done in the \rdga, a set of candidate orders is maintained, from which the next order is selected. After each selection, new candidate orders are identified, and the set of candidate orders is re-sampled from the combined old and new orders. New candidate orders are picked based on custom distance functions based on the newly considered articles belonging to the added order.

In our testing implementation, we considered a base set of $100$ candidate orders, with $10$ new candidate orders being considered in each step. In each step, from these candidate orders, the one closest to already picked items is selected. Then $10$ new orders are identified whose items are in the same aisle as previously picked items. This aims to achieve synergies, as some of these items might lie on existing routes already and if not, extending these routes to such items is cheap on average.

\subsection{Computational Evaluation - Single Phase}\label{sec:comp-single}
Before we present the performance evaluation data, we give a brief overview of the instance parameters and instance generation, for additional details we refer to Appendix \ref{sec:app:instanceparams}.
We provide additional data for other instances in Appendix \ref{sec:app:singlephasedata}.

Each instance is created in a way that ensures that all orders can always be fulfilled. Order sizes and article volumes are sampled from appropriate probability distributions and each order consists of randomly chosen articles. The number of items in the warehouse is roughly three times as large as the number of articles and each item is assigned a uniformly random location in the warehouse. We present data for two different instance sizes called \emph{medium} and \emph{large} with parameters outlined in \cref{tab:instanceparams}.
These instance parameters were also taken from \cite{ZalandoPaper} and later generalized for our multi-phase approach.
For each instance class, the item goal is set to be approximately 20\% of the total number of articles on order.

\begin{table}
  \TABLE
  {Parameters for the different instance categories. \label{tab:instanceparams}
  }
  {
    \begin{tabular}{l@{\hskip 2\tabcolsep}|c@{\hskip 6\tabcolsep}c@{\hskip 6 \tabcolsep}c@{\hskip 6\tabcolsep}}\toprule
      Parameter &  medium & large \\\midrule
      Warehouse items &  100,000 & 1,000,000\\
      Orders &  5,000 & 50,000\\
      Zones & 50 & 100\\
      Rows &  100 & 100\\
      Aisles & 100 & 100 \\\bottomrule
  \end{tabular}}
  {}
\end{table}

We introduce the following variations on some parameters of the instance descriptions: by ``\nobreakdash-fewzones'' and ``-fewerzones'', we denote the instance that uses only 10 and 3 zones, respectively.
Another variant is the ``-longbatch'' instance, in which 250 orders may be scheduled per batch, as compared to the 50 orders in the standard instances.

\subsubsection{Analysis of the computational results}
\label{sec:singlephase-empirical-analysis}
In this section we computationally evaluate all algorithms outlined above and analyze their advantages and disadvantages. For \rdgak, we choose two variants: $k=50$ and $k= \sqrt{|\Orders|}$.
Each algorithm was applied to the same set of five randomly generated instances of medium and large size. The numerical data indicates the average value of the respective observation over all runs of the same instance size.
As the variance between runs on the same instance and over the different instances of the same size is low, no further statistical parameters other than the mean value are presented.

In \cref{tab:singlephaseperformance}, the results for the primary objectives in our single-phase analysis are displayed: the total cost of the picklists and the runtime of their computation.
Note that the cost of \dga and \fdga is not only equal on average, but their output is exactly the same, and we therefore only provide data for \fdga outside of runtime discussions.
Our first contribution is the vastly improved runtime, which is about $18$ times faster for \fdga compared to \dga in large instances, thereby making this approach significantly more feasible in terms of real-world applicability.
We observe that our first variant \fdgaskip is able to reduce the runtime by a further 8--10\% with a total reduction of a factor of $19.5$ compared to \dga, while maintaining solution quality that is comparable to the reference algorithm (within $\pm 0.3\%$ of the objective for large instances).

The \rdga variants implement a sampling approach to trade some of the cost performance in favor of even faster runtime, which we observe in this data: While the increase in cost compared to the reference is not insignificant, the computation time is another order of magnitude faster than \fdga.
The \ldga algorithm offers a compromise between the two approaches, with increased runtime compared to \rdga due to the overhead incurred for maintaining candidate orders, but at the same time achieving cost closer to \fdga (+14\%) than the independent sampling approaches.

\begin{table}
  \TABLE{Benchmark data of selected algorithms. \label{tab:singlephaseperformance}}
  {\footnotesize
    \begin{tabular}{c@{\hskip 2\tabcolsep}|ccr|ccr}\toprule
      & \multicolumn{3}{c|}{medium} & \multicolumn{3}{c}{large} \\
      Algorithm  & Cost                        & Relative Performance & Time & Cost                        & Relative Performance & Time \\\midrule
      \dga & 62354.4 & +0.2\% & 108.9s & 524925.2 & +0.0\% & 10548.6s \\
      \fdga & 62354.4 & +0.2\% & 8.7s & 524925.2 & (reference) & 591.8s \\
      \fdgaskip & 62203.2 & (reference) & 7.9s & 526366.8 & +0.3\% & 545.5s \\
      \ldga & 70684.4 & +13.6\% & 5.7s & 597082.8 & +13.7\% & 113.8s\\
      \rdgasqrt & 76453.2 & +22.9\% & 1.8s & 734645.2 & +40.0\% & 59.1s \\
      \rdgan{50} & 77449.2 & +24.5\% & 1.4s & 843182.0 & +60.6\% & 22.5s\\
      \rdga & 139287.2 & +123.9\% & 0.2s& 1489700.8 & +183.8\% & 9.8s \\\bottomrule
  \end{tabular}}{All values are averaged over five different instances and rounded to the nearest tenth.}
\end{table}

\subsection{Further Solution Insights}
\label{sec:singlephase-performance}

In this section we analyze additional solution properties. The goal is to gain insight into the solution structure and we perform additional experiments to evaluate different properties of the optimization procedure.

\begin{table}
  \TABLE{Further computational measures for different algorithms.
  \label{tab:singlephaseinsights}}
  {\footnotesize
    \begin{tabular}{c@{\hskip 2\tabcolsep}|ccc|ccc}\toprule
      & \multicolumn{3}{c|}{medium} & \multicolumn{3}{c}{large} \\
      Algorithm & \shortstack{Item distance\\from conveyor}  & Order size & Picklist length & \shortstack{Item distance\\from conveyor}  & Order size & Picklist length \\\midrule
      \fdga & 23.2 & 2.3 &  3.3 & 22.8 & 2.2 &  2.6 \\
      \fdgaskip & 23.3 & 2.3 & 3.3 & 22.8 & 2.2 &  2.6 \\
      \ldga & 23.3 & 2.3 & 3.3 & 23.2 & 2.3 &  2.5 \\
      \rdgasqrt& 22.1 & 2.4 &  3.0 & 21.5 & 2.4  & 2.1 \\
      \rdgan{50}& 22.0 & 2.4 &  2.9 & 21.1 & 2.4 & 2.0 \\
      \rdga & 33.5 & 2.7 &  3.0& 32.8 & 2.6  & 1.9 \\
      (Instance) & 34.03 & 2.65 & $-$ &  34.03 & 2.64  &$-$    \\\bottomrule
  \end{tabular}}{All values are averaged over five different instances and represent the mean value over all batched items, picklists, and orders respectively. The reference values refer to all orders/items in the instance as compared to the ones picked by the algorithm.}
\end{table}

\cref{tab:singlephaseinsights} outlines some observations for the previously discussed algorithms.
The picklist length is the average length over all picklists created in the last step of each algorithm.
The order size is the average number of items contained in the selected orders, and the item distance from the conveyor is the average distance each of the picked items has to the conveyor.
The instance reference value for the average item distance from the conveyor belt refers to the resulting distances when all available orders are fulfilled.
This reference has been computed by selecting the closest of the available items of each article to fulfill each order.
Because we do not use actual route lengths for this estimate, the specific item selection does not affect this metric.

We observe that the average distance to the conveyor of items picked by the algorithms is significantly lower than the average over all items (with the exception of the crude random item choice in \rdga). While some of this difference is to be expected due to the algorithm usually selecting the item corresponding to an article which is closest to the conveyor, the observed effect is actually stronger than this. It turns out that the algorithms strongly prefer orders whose items are close to the conveyor, which when considering the multi-phase setting leads to a kind of gradual emptying of the center of each zone and a simultaneous starvation of orders whose items are mostly in far positions.
Additionally, the length of batched orders is less than the average, indicating that the algorithms favor batching smaller orders, as it might be easier to accommodate all items in an order if there are fewer of them.
This behavior is obviously undesirable from a managerial perspective, as it could lead to unbalanced warehouse utilization and delayed order fulfillment. We capture and address these observations in our extension to multiple phases discussed in \cref{sec:multiphase}.

\paragraph{Picklist Length and Routing}
Another observation is that the number of items picked in each picklist is quite small, even for the reference algorithm, with on average fewer than 4 items picked for each picker route.
This is due to the restrictive batch size (50 orders) and small average order size in comparison to the number of zones (50 for medium, 100 for large instances) in the problem.

To elaborate on this observation and eliminate this undesirable effect, we have performed additional computations on varied parameters as outlined above and provide some performance-related data in \cref{tab:singlephaseperformance-fewerzones}.
The performance evaluation confirms the qualitative ranking of the different algorithms outlined above.
Additionally, we observe a more realistic outcome with longer picklists and higher load on the pickers' capacity constraints: for the -fewzones variants, the picklists contain 10--13 items on average with a median picker load of ca. 30--40\%, and for the -fewerzones variants, the picklists contain 17.5--20 items per picklist with median load exceeding {80\%}.

\begin{table}
  \TABLE{Benchmark data of selected algorithms for medium instances with few and fewer zones. \label{tab:singlephaseperformance-fewerzones}}
  {
    \begin{tabular}{c@{\hskip 2\tabcolsep}|ccr|ccc}\toprule
      & \multicolumn{3}{c|}{medium-fewzones} & \multicolumn{3}{c}{medium-fewerzones} \\
      Algorithm  & Cost                        & Relative Performance & Time & Cost                        & Relative Performance & Time \\\midrule
      \dga & 41446.4 & +0.0\% & 185.5s & 29281.6 & +2.7\% & 395.1s \\
      \fdga & 41446.4 & (reference) & 30.4s & 29281.6 & +2.7\% & 74.7s \\
      \fdgaskip & 41682.4 & +0.6\% & 27.7s& 28498.4 & (reference) & 67.7s \\
      \ldga & 43815.6 & +5.7\% & 9.1s& 33404.4 & +17.2\% & 21.1s \\
      \rdgan{50} & 53383.2 & +28.8\% & 2.3s& 38494.0 & +35.1\% & 4.7s \\
      \rdgasqrt & 51320.4 & +23.8\% & 2.9s & 36911.6 & +29.5\% & 6.1s \\
      \rdga & 102917.6 & +148.3\% & 0.2s& 83658.4 & +193.6\% & 0.3s \\\bottomrule
  \end{tabular}}{ All values are averaged over five different instances and rounded to the nearest tenth.}
\end{table}

While \cite{ZalandoPaper} have only used a heuristic for the computation of the picklists after the item selection has concluded, we have additionally experimented with improving solutions for this subproblem using a specialized library \citep{ortools_routing} for constrained vehicle routing problems (CVRP) for the optimization separately on every zone.
For this solver, we have used the default parameters with a time limit of 10s for each zone.

Even though the basic routing heuristic performs sufficiently well in the reference instances due to the small number of items in each picklist and only provides marginal improvements in the objective value (\cref{tab:singlephasecvrp-comparison}), we observe a different outcome with substantially improved performance in our problem variant with longer picklists.
For later results, we therefore restrict our analysis to routing using the external library instead.

\begin{table}
  \TABLE{Comparing the heuristic routing subroutine with the improved CVRP-based solution for medium and medium-fewerzones instances.
  \label{tab:singlephasecvrp-comparison}}
  {
    \begin{tabular}{c@{\hskip 2\tabcolsep}|ccc|ccc}      \toprule

      & \multicolumn{3}{c|}{medium} & \multicolumn{3}{c}{medium-fewerzones} \\

      Algorithm  & Greedy    & CVRP  & Rel. Performance & Greedy    & CVRP  & Rel. Performance \\
      \midrule
      \fdga & 62354.4  & 61526.4  & -1.3\% & 29281.6 & 26103.2  & -10.9\% \\
      \fdgaskip & 62203.2  & 61215.6  & -1.6\% & 28498.4  & 25474.0  & -10.6\% \\
      \ldga & 70684.4  & 68294.4  & -3.4\% & 33404.4  & 29532.8  & -11.6\% \\
      \rdgasqrt & 76453.2  & 75569.2 & -1.2\% & 36911.6  & 34086.0  & -7.7\% \\
      \rdgan{50} & 77449.2  & 77311.6  & -0.2\% & 38494.0  & 35381.2  & -8.1\% \\
      \rdga & 139287.2  & 135510.8  & -2.7\% & 83658.4  & 64500.4  & -22.9\% \\
      \bottomrule
  \end{tabular}}{All values are averaged over five different instances and rounded to a tenth.}
\end{table}

\section{Multi Phase Setting}
\label{sec:multiphase}
In this section we discuss and analyze the extension of the model to multiple phases, as defined in \cref{sec:multiphase-model}. This extension is motivated by the need to more accurately account for settings observed in real-life warehouses and to thus reduce weaknesses of the single-phase model.

\subsection{Motivation}
Each picking phase in the original single-phase model is supposed to represent a rather short timeframe during a typical workday at a warehouse. Depending on the exact circumstances, the represented timeframe might range from a few minutes up to around an hour. So a typical workday will consist of many such phases, with restockings of the warehouse throughout. While restocking might actually be performed in parallel to order picking in some warehouses, we slightly simplify the model and assume instantaneous restocking in-between picking phases. Due to the rather short length of these phases, this does not have a large impact.

Considering multiple consecutive phases has the advantage of both capturing possible synergies between picking decisions over multiple phases while also eliminating a biased advantage which single-phase algorithms can obtain by making short-term decisions. Some such advantages were already observed in \cref{sec:singlephase-performance}, as many of the algorithms overly prioritize orders which can be fulfilled by items close to the conveyor belt or which only consist of few articles.
We show that over a longer timeframe, such decisions lead to worse overall outcomes, as later phases are then much harder to fulfill cost efficiently. The multi-phase setting is also beneficial from a managerial perspective, as other performance indicators such as fulfillment guarantees and order starvation prevention can be tracked, analyzed, and optimized for.

As a motivating example, we would like to present a glimpse of the data outlined in \cref{sec:multiphase-empirical-analysis}: After the naive application of \fdga to 100 consecutive phases on medium instances, we end up with an order set, where the size of an order yet to be fulfilled is 3.6 articles on average (compared to {2.64} articles of average order size in the ground set) and the average age of an order (in terms of number of phases since it has been added to the problem) is 15.72 phases. The oldest order that has not been fulfilled yet is from the initial phase, i.e.\ it has been added 100 phases prior.
Under this light, the competitive performance of \fdga must be critically evaluated, as it can be conjectured that it achieves smaller cost only at the trade-off of leaving generally unsynergetic orders unfulfilled.

In the following, we present algorithms which not only eliminate this potential drawback, but in some cases even outperform the benchmark \fdga algorithm in terms of the total cost. Furthermore, this is all achieved without increasing the runtime of \fdga.

\subsection{Computational Evaluation - Multi Phase}
Similar to the single-phase setting, we evaluate the different algorithms regarding their cost and runtime. In this case, we consider $100$ consecutive phases. Each phase is treated by the algorithms as an independent instance of a single-phase picking problem. The only additional information that is communicated to the solvers is the \emph{age} of every order currently in the warehouse, i.e.\ the number of phases that have passed since the order came in.

In order to computationally model the sequential setting, we devised a scheme to populate the warehouse with items and orders similar to the single-phase case. Initially, we start with some number of orders that vastly exceeds the number of orders to be completed when the item goal of a single phase is approximately fulfilled, and a warehouse stocked with some number of extra items that are currently not requested, just as in the single-phase setting.
On a technical level, we generate orders and their according items as sequences. After the completion of each phase, we take the next chunk of items and corresponding orders that they can fulfill from the sequence and add them to the instance. The full sequences are unknown to the algorithm and new items and orders are only revealed during the restocking.
In order to maintain a consistent warehouse size, we restock the exact number of items that was removed in the previous phase. These items also approximately correspond to the number of orders completed in the previous phase.
\Cref{tab:instanceparams} outlines the parameters used in this experiment, and they are directly applied to the initial phase of the multi-phase setting and thus generalize the instance categories \emph{medium} and \emph{large} to this setting. Within these parameters, approximately $20\%$ of the available orders are fulfilled in every phase, i.e.\ around 1,000 orders for the medium instances and 10,000 for the large instances.

\subsection{Multi-Phase Algorithms}

While the single-phase setting did not have a notion of \emph{age} for each order, this is a feature we consider crucial in our extension to the multi-phase setting. As outlined below, this is also reflected in our newly devised algorithms. Additionally, we compare the novel variations to the iterative and age-agnostic application of selected algorithms (\fdga, \ldga, and \rdgasqrt) from the single-phase setting above, which compute independent instances for each new phase. As \fdga computes the same solutions as \dga, this selection thus provides a benchmark to the application of the algorithm from \cite{ZalandoPaper} while drastically improving the computation time.

\subsubsection{\fdga with soft age objective}

For this first approach of incorporating the age information of each order to achieve timely order fulfillment, we modify the cost function in the \fdga procedure to prioritize old orders over new ones to a certain degree.
In order to do this, we discount the cost of older orders according to their age using two parameters $\gamma$ and $\alpha$, where the former corresponds to geometric and the latter to linear discounts.
The cost heuristic in \fdga which selects the next order now computes the cost $c(O)$ for each order as before and modifies it to a new age-dependent cost of $\tilde c(O) = \gamma^t c(O) - \alpha t$, where $t$ is the age of order $O$ given as an integer counting the number of phases.
The order that minimizes $\tilde c$ is now returned.
Thus, all orders of the same age are again compared to each other with the same discount factor.
Depending on the given parameter choice of $\gamma$ and $\alpha$, older orders are prioritized over new ones, while still maintaining the incentive for the algorithm to choose synergising orders.

In our experiments, we have found satisfactory performance without using any linear discounts (i.e.\ $\alpha=0$) and refer to the algorithms by \fdgasoftg{$\gamma$} for the corresponding geometric discount parameter $\gamma$.
We use $\gamma \in \{0.8,\ 0.9,\allowbreak 0.95,\ 0.99\}$ for our experiments, which corresponds to discounting the cost of each order for each phase by $20\%$, $10\%$, $5\%$, and $1\%$ respectively.

\subsubsection{\fdga with strict age objective}

While the previous approach accomplishes the completion of old orders implicitly due to the modification of the costs, we have also devised a second approach that gives explicit hard guarantees by design of the algorithm. For a full description see Appendix \ref{sec:app:fdgastrictage}.

This variant of \fdga takes as input a list of non-increasing order numbers $B = b_0,\dots, b_m$, where $b_m = 0$.
We also call the $b_i$ \emph{bucket size of age $i$} or \emph{completion targets}.
Each $b_i$ indicates how many orders of age $i$ may be carried over to the next phase of the algorithm.
Thus, no order can exceed age $m$ over the course of the algorithm. For appropriately chosen bucket sizes, the algorithm is guaranteed to hit the target age distributions. This is achieved by randomly sampling as many orders as necessary to fulfill the completion targets and then employing a modified version of \fdga to select both these mandatory orders as well as additional synergistic orders in an efficient way.

In our experiments, we use the following implementations of the bucket sizes. Note that both of them implicitly use (an estimate of) the number $n$ of newly added orders in each phase.

\begin{itemize}
  \item \fdgalinearn{$m$} implements a linear decline of bucket sizes from $n$ to 0 at age $m$.
    The buckets are thus given by $b_i = \max(0, \lceil n\cdot (m-i)/m\rceil)$
  \item \fdgacutoffn{$m$} does not impose any restrictions on the orders up to age $m$, but then requires all orders of that age to be completed.
    Thus we have $b_0,\dots, b_{m-1} = n$ and $b_m = 0$.
\end{itemize}

\subsection{Multi-Phase Analysis}
\label{sec:multiphase-empirical-analysis}

First, we start off by comparing the performance of the different algorithms in \cref{tab:multiphase-performance-medium-fewzones} for medium-fewzones instances, as we deem its solutions an adequate compromise between the number of zones to pose a sufficiently challenging problem and the resulting picklist lengths for practical applications, as discussed in \cref{sec:comp-single}.
We provide analogous data for different instance variants in Appendix \ref{sec:app:multiphasedata}.
All of the results reflect our qualitative observations outlined in this section, but the quantitative extent varies slightly between the specific instance settings.

Besides the total cost, we also present the relative cost of these algorithms compared to the single-phase reference \fdga.
Additionally, the average cost per batched item and per batched order is given, in order to eliminate any irregularities from varying numbers of chosen items or orders due to the order selection in the computation.

For different parameter choices for the algorithms \fdgacutoffn{$m$} and \fdgalinearn{$m$} variants we can clearly observe the decreased opportunities for synergy during the computation, when the completion targets are set more tightly.
Nevertheless, we do observe only marginally increased cost for the trade-off with strict completion guarantees, specifically when compared to the application of the more runtime efficient single-phase algorithm \rdgasqrt.
On the contrary, we even observe \emph{decreased} cost for most \fdgasoft parameter choices compared to the benchmark \fdga.

\begin{table}
  \TABLE{Benchmark data comparing different multi-phase algorithms for medium-fewzones instances.
    \label{tab:multiphase-performance-medium-fewzones}
  }{\footnotesize
    \setlength{\tabcolsep}{3pt}
    \begin{tabular}{c@{\hskip 2\tabcolsep}|cccc|cccc}\toprule
      Algorithm & Total cost & Rel.\ cost & \shortstack{Avg.\ cost\\ per item} & \shortstack{Avg.\ cost\\ per order} & \shortstack{Average\\ order length} & \shortstack{Maximum\\ order age} & \shortstack{Average\\ order age} & \shortstack{Median\\ order age} \\\midrule
      \fdga & 5609744 & (reference) & 21.51 & 55.97 & 3.6 & 100.0 & 15.72 & 4.2 \\
      \ldga & 5566124 & $-0.8\% $ & 21.41 & 55.9 & 3.31 & 100.0 & 12.35 & 3.0 \\
      \rdgasqrt & 6620678 & $+18.0\% $ & 25.44 & 66.5 & 3.21 & 100.0 & 12.11 & 4.0 \\
      \fdgacutoffn{5} & 7023362 & $+25.2\% $ & 27.28 & 71.84 & 2.73 & 5.0 & 2.34 & 2.0 \\
      \fdgacutoffn{10} & 6652562 & $+18.6\% $ & 25.75 & 67.46 & 3.03 & 10.0 & 3.53 & 3.0 \\
      \fdgacutoffn{15} & 6352230 & $+13.2\% $ & 24.59 & 64.29 & 3.17 & 15.0 & 4.56 & 3.0  \\
      \fdgacutoffn{20} & 6188696 & $+10.3\% $ & 23.87 & 62.34 & 3.26 & 20.0 & 5.46 & 3.0\\
      \fdgacutoffn{25} & 6067686 & $+8.2\% $ & 23.4 & 61.07 & 3.32 & 25.0 & 6.27 & 3.2\\
      \fdgacutoffn{30} & 5991385 & $+6.8\% $ & 23.04 & 60.12 & 3.38 & 30.0 & 7.05 & 3.4 \\
      \fdgalinearn{10} & 6750239 & $+20.3\% $ & 26.15 & 68.55& 3.0 & 10.0 & 3.37 & 3.0 \\
      \fdgalinearn{15} & 6423991 & $+14.5\% $ & 24.81 & 64.88 & 3.16 & 15.0 & 4.47 & 3.0\\
      \fdgalinearn{20} & 6226611 & $+11.0\% $ & 24.03 & 62.77 & 3.26 & 20.0 & 5.43 & 3.0\\
      \fdgalinearn{25} & 6113034 & $+9.0\% $ & 23.53 & 61.43 & 3.33 & 25.0 & 6.31 & 3.2\\
      \fdgalinearn{30} & 6018649 & $+7.3\% $ & 23.16 & 60.43 & 3.38 & 30.0 & 7.11 & 3.8\\
      \fdgasoftg{0.8} & 5539346 & $-1.3\% $ & 21.39 & 56.25 & 2.81 & 8.2 & 2.29 & 2.0\\
      \fdgasoftg{0.9} & 5315791 & $-5.2\% $ & 20.55 & 53.93& 2.93 & 12.6 & 2.71 & 2.0  \\
      \fdgasoftg{0.95} & 5299643 & $-5.5\% $ & 20.43 & 53.51& 3.06 & 20.0 & 3.38 & 2.0 \\
      \fdgasoftg{0.99} & 5469906 & $-2.5\% $ & 21.01 & 54.84 & 3.35 & 60.6 & 6.49 & 3.0\\\bottomrule
  \end{tabular}}{The parameter for the number of new orders is set to 1000 if applicable and approximately corresponds to the number of new orders in each phase. The column for relative cost compares to the reference of \fdga. The values are averaged over 5 instances and rounded to appropriate precision. The values for the second part of the table are computed for the set of orders in the \emph{last} phase of each run. For reference, the average order length in the instance is 2.64.}
\end{table}

\begin{figure}
  \FIGURE{
  \includegraphics[width=0.8\textwidth]{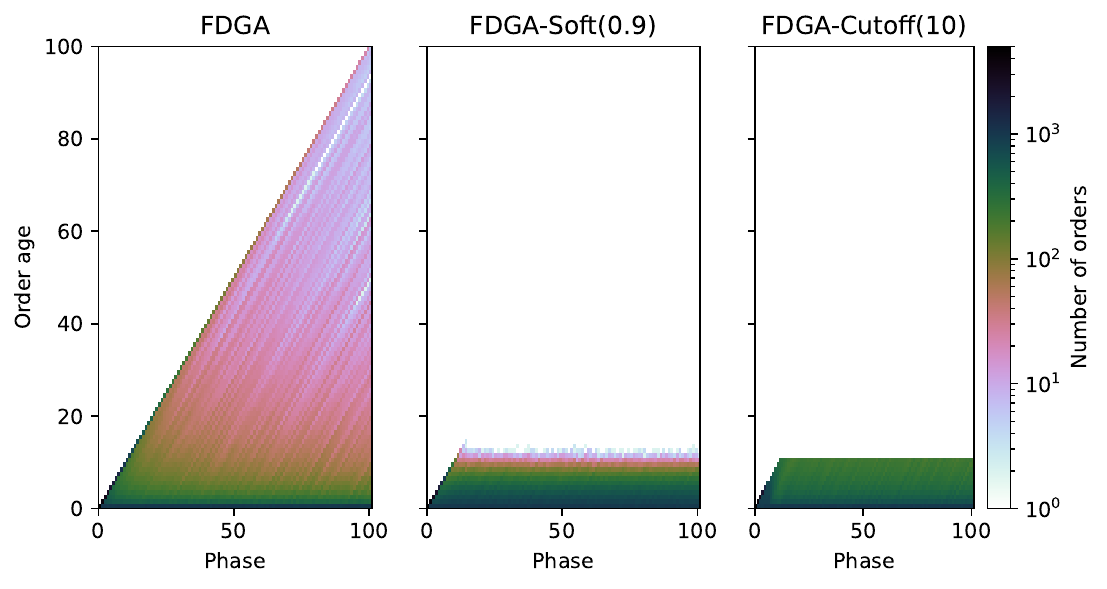}}
  {Age distribution data of single runs of \fdga, \fdgasoftg{0.9}, and \fdgacutoffn{10} on a medium-fewzones instance.
  \label{fig:fdga_age}}{While many orders are completed right away, many of them remain in the problem until the end of the experiment.}
\end{figure}

Next, we compare the solutions more thoroughly than just with regards to the cost performance, as our original goal was in fact not to improve the performance over the benchmark, but rather to produce \emph{better} solutions in terms of other properties.
A crucial observation on the solution structure of the algorithms is illustrated in \cref{fig:fdga_age}.
The application of the standard \fdga algorithm completes the majority of newly added orders when they are added to the problem, but the set of remaining orders is not necessarily completed in reasonable time.
More specifically, we observe pronounced diagonal lines which barely change over time, which represent the orders added in a given phase that are apparently ignored throughout the rest of the phases and increase the size of the backlog, as the algorithm is not able to accommodate them according to its objectives.
This behavior motivated the development of our novel algorithms which fix this issue by design.

In \cref{fig:fdga_age}, we also provide exemplary data for the completion ages of orders for \fdgasoft and \fdgacutoff, where all orders are being completed within around 10 phases after being added to the problem.
Comparing these two algorithms against each other, we can also confirm that they accomplish the design goals: While \fdgasoft prioritizes older orders over newer ones and their number declines systematically as they age, the distribution of ages remains more consistent in \fdgacutoff until all the respective orders are mandatorily completed at the cutoff age.

The same qualitative behavior can be observed for other parameter choices in these algorithms and we discuss them numerically next.
Further data for this experiment in the multi-phase setting is displayed in the later columns in \cref{tab:multiphase-performance-medium-fewzones}.
This data is computed from the set of orders within the last phase of the multi-phase model, again averaged over all 5 independent runs.
As such, we believe that it is representative of the convergence behavior of the application of these algorithms over time (where we suggest that some of the discussed effects could be even \emph{more} drastic for longer time periods).

With regards to the current state of the warehouse in terms of orders, we list the average number of positions in each order in the warehouse.
This confirms a similar observation from the single-phase experiments, where some algorithms strongly prefer shorter orders over longer orders, most likely due to the fact that fewer items are more likely to fit in well with the already selected routes.
This is reflected very well in the comparison across different parameters for our new algorithm variants, where for stronger completion target choices, the proportion of longer orders in the warehouse is closer to the instance average.

Another interesting observation is the distribution of ages of orders, where we continue to observe very regular behavior in the parameters for the different variants.
We can empirically confirm that the completion targets for \fdgalinear and \fdgacutoff work as intended, while we also get a rough estimate for the effective upper bound on the targets for \fdgasoft from the data.
Most notably, we want to emphasize that the problematic behavior of non-fulfilled orders from the initial phase is confirmed by our experiments.
The average order age is well controlled only by our newly introduced variants, while the effect is observable, but less pronounced in the median order age.

Another observation that we do not present numerically is that all of the observed algorithms prefer to pick items closer to the center of the warehouse.
Due to the design of the objective, this is also expected behavior and we did not observe any notable differences between the algorithms in this regard.

\section{Conclusion}
Our contribution builds upon the model and algorithms \dga and \rdga from \cite{ZalandoPaper}, which we reproduce in \cref{sec:singlephase-model}.
For our own implementation of \fdga, we improve upon \dga by efficiently using both data structures and problem-specific properties to significantly reduce the runtime.
By skipping the computation of a subroutine, we observe qualitatively comparable outcomes with even faster runtime in the algorithm \fdgaskip.
Furthermore, we introduce a variant called \rdgak, which allows the user to compromise between runtime and performance via a single parameter, while experimentally maintaining decent solution quality compared to the reference of \dga for reasonable fixed parameter choices or a dynamic and instance dependent approach with \rdgasqrt.
Additionally, we present novel algorithmic approaches aiming to provide different (and configurable) trade-offs between runtime and cost.
Thorough experimental evaluation of all approaches is contained in \cref{sec:singlephase}.

We further extend the model by \cite{ZalandoPaper} to a multi-phase setting, more closely resembling the reality observed in warehouse operations. We show that a naive extension of the single-phase algorithms to this setting returns solutions with undesirable characteristics, such as arbitrarily long order delays and unbalanced warehouse utilization.
Our novel algorithms address these shortcomings and are presented and analyzed in \cref{sec:multiphase}.
We show that surprisingly simple modifications can lead to solutions with strong guarantees, while maintaining solution quality or even improving upon the reference.

For both the single phase as well as the multi-phase model we additionally analyze instance variants which represent denser demand profiles with higher picker utilization than the ones given in \cite{ZalandoPaper}. We show that our results extend to this setting and are generally robust against input variations.

A limitation in this work is the simplified restocking model employed in the multi-phase setting.
In our model, the restocked items are placed uniformly at random throughout the warehouse and are drawn uniformly from the total pool of available items.
In real world settings, additional data can provide estimates for both the items being re-stocked as well as the items being ordered, which is taken into account during warehouse design and when deciding where to place newly arrived items. Incorporating such decisions and procedures into the algorithms could lead to further overall improvements in terms of costs and fulfillment guarantees. In a similar vein, extensions of the analysis performed in this work to other warehouse settings and structures as discussed in \cref{sec:intro} are of interest. Finally, testing both the algorithms from this work as well as potential variants for other settings using real historic data from warehouse operations would be an important step towards real world deployment.

\bibliographystyle{alpha}
\bibliography{bib}

\newcommand{\etalchar}[1]{$^{#1}$}
\begin{thebibliography}{GVdBVDH01}

\bibitem[AA{\"O}25]{Aksu2025-ordersize5}
Ozan~R{\i}dvan Aksu, Necati Aras, and Temel {\"O}ncan.
\newblock Solving the joint order batching and picker routeing selection
  problem in warehouses.
\newblock {\em Int. J. Prod. Res.}, pages 1--26, October 2025.

\bibitem[ACS21]{Aerts2021-routeassignment}
Babiche Aerts, Trijntje Cornelissens, and Kenneth S{\"o}rensen.
\newblock The joint order batching and picker routing problem: Modelled and
  solved as a clustered vehicle routing problem.
\newblock {\em Comput. Oper. Res.}, 129(105168):105168, May 2021.

\bibitem[AKK{\etalchar{+}}25]{ZalandoPaper}
Giorgio Abelli, Maximilian Katzmann, Imran Khan, Olaf Maurer, Julius Pätzold,
  Paweł Pszona, and Jan-David Salchow.
\newblock Joint order selection, allocation, batching and picking for large
  scale warehouses, 2025.

\bibitem[BBBS24]{Bock2024-WarehouseSurvey}
Stefan Bock, Stefan Bomsdorf, Nils Boysen, and Michael Schneider.
\newblock A survey on the traveling salesman problem and its variants in a
  warehousing context.
\newblock {\em Eur. J. Oper. Res.}, April 2024.

\bibitem[BBV{\etalchar{+}}24]{Bocelli2024-warehousedesign}
Michele Bocelli, Eleonora Bottani, Andrea Volpi, Federico Solari, Natalya
  Lysova, and Roberto Montanari.
\newblock Exploring the relationship between routing policies and market demand
  heterogeneity.
\newblock {\em J. Ind. Inf. Integr.}, 40(100634):100634, July 2024.

\bibitem[BCC{\etalchar{+}}20]{Briant2020-routeassignment}
Olivier Briant, Hadrien Cambazard, Diego Cattaruzza, Nicolas Catusse,
  Anne-Laure Ladier, and Maxime Ogier.
\newblock An efficient and general approach for the joint order batching and
  picker routing problem.
\newblock {\em Eur. J. Oper. Res.}, 285(2):497--512, September 2020.

\bibitem[BdKW19]{Boysen2019-WarehousingInECommerce}
Nils Boysen, Ren{\'e} de~Koster, and Felix Weidinger.
\newblock Warehousing in the e-commerce era: A survey.
\newblock {\em Eur. J. Oper. Res.}, 277(2):396--411, September 2019.

\bibitem[CC18]{Cambazard2018-np}
Hadrien Cambazard and Nicolas Catusse.
\newblock Fixed-parameter algorithms for rectilinear steiner tree and
  rectilinear traveling salesman problem in the plane.
\newblock {\em Eur. J. Oper. Res.}, 270(2):419--429, October 2018.

\bibitem[CCCCE21]{Cano2021-ordersize3}
Jose~Alejandro Cano, Pablo Cort{\'e}s, Emiro~Antonio Campo, and
  Alexander~Alberto Correa-Espinal.
\newblock Solving the order batching and sequencing problem with multiple
  pickers: A grouped genetic algorithm.
\newblock {\em Int. J. Electr. Comput. Eng. (IJECE)}, 11(3):2516, June 2021.

\bibitem[CLSG24]{Cui2024-deliveryadvantage}
Ruomeng Cui, Zhikun Lu, Tianshu Sun, and Joseph~M Golden.
\newblock Sooner or later? promising delivery speed in online retail.
\newblock {\em Manuf. Serv. Oper. Manag.}, 26(1):233--251, January 2024.

\bibitem[CSRM15]{Cardona2015-warehousedesign}
Luis~F Cardona, Diego~F Soto, Leonardo Rivera, and Hector~J Mart{\'\i}nez.
\newblock Detailed design of fishbone warehouse layouts with vertical travel.
\newblock {\em Int. J. Prod. Econ.}, 170:825--837, December 2015.

\bibitem[CVM{\etalchar{+}}23]{CVMTB:TrendsInOrderPicking}
Giorgia Casella, Andrea Volpi, Roberto Montanari, Letizia Tebaldi, and Eleonora
  Bottani.
\newblock Trends in order picking: a 2007–2022 review of the literature.
\newblock {\em Production \& Manufacturing Research}, 11(1), 2023.

\bibitem[CZLZ23]{Cao2023-ordersize2}
Zhengcai Cao, Lijie Zhou, Chengran Lin, and Mengchu Zhou.
\newblock Solving an order batching, picker assignment, batch sequencing and
  picker routing problem via information integration.
\newblock {\em J. Ind. Inf. Integr.}, 31(100414):100414, February 2023.

\bibitem[DSL{\etalchar{+}}25]{dunn2025det}
George Dunn, Elizabeth Stojanovski, Bishnu Lamichhane, Hadi Charkhgard, and Ali
  Eshragh.
\newblock Deterministic structure of vertical configurations in minimal picker
  tours for rectangular warehouses, 2025.

\bibitem[FGX19]{Fisher2019-deliveryadvantage}
Marshall~L Fisher, Santiago Gallino, and Joseph~Jiaqi Xu.
\newblock The value of rapid delivery in omnichannel retailing.
\newblock {\em J. Mark. Res.}, 56(5):732--748, October 2019.

\bibitem[FP26]{ortools_routing}
Vincent Furnon and Laurent Perron.
\newblock Or-tools routing library, 2026.

\bibitem[GV05]{Gademann2005-jg}
Noud Gademann and Steef Velde.
\newblock Order batching to minimize total travel time in a parallel aisle
  warehouse.
\newblock {\em IIE Trans.}, 37(1):63--75, January 2005.

\bibitem[GVdBVDH01]{Gademann2001-cv}
Noud Gademann, Jeroen~P Van~den Berg, and Hassan~H Van Der~Hoff.
\newblock An order batching algorithm for wave picking in a parallel-aisle
  warehouse.
\newblock {\em IIE Trans.}, 33(5):385--398, May 2001.

\bibitem[Kar72]{Karp1972-21probs}
Richard~M Karp.
\newblock Reducibility among combinatorial problems.
\newblock In {\em Complexity of Computer Computations}, pages 85--103. Springer
  US, Boston, MA, 1972.

\bibitem[KGB20]{Kubler2020-ordersize1}
Patrick K{\"u}bler, Christoph~H Glock, and Thomas Bauernhansl.
\newblock A new iterative method for solving the joint dynamic storage location
  assignment, order batching and picker routing problem in manual
  picker-to-parts warehouses.
\newblock {\em Comput. Ind. Eng.}, 147(106645):106645, September 2020.

\bibitem[LBS22]{Loffler2022-AGVPicking}
Maximilian L{\"o}ffler, Nils Boysen, and Michael Schneider.
\newblock Picker routing in {AGV-assisted} order picking systems.
\newblock {\em INFORMS J. Comput.}, 34(1):440--462, January 2022.

\bibitem[MGG20]{Masae2020-OrderPickingLitRev}
Makusee Masae, Christoph~H. Glock, and Eric~H. Grosse.
\newblock Order picker routing in warehouses: A systematic literature review.
\newblock {\em Int. J. Prod. Econ.}, 224(107564):107564, June 2020.

\bibitem[MGV21]{Masae2021-ordersizeleaf}
Makusee Masae, Christoph~H Glock, and Panupong Vichitkunakorn.
\newblock A method for efficiently routing order pickers in the leaf warehouse.
\newblock {\em Int. J. Prod. Econ.}, 234(108069):108069, April 2021.

\bibitem[{{Oxf}}24]{ZalandoReport2024}
{{Oxford Economics}}.
\newblock The economic impact of zalando in europe.
\newblock Technical report, Oxford Economics, Zalando SE, Marienstraße 15,
  60329, Frankfurt, Germany, December 2024.

\bibitem[OZD17]{Onal2017-explosive}
Sevilay Onal, Jingran Zhang, and Sanchoy Das.
\newblock Modelling and performance evaluation of explosive storage policies in
  internet fulfilment warehouses.
\newblock {\em Int. J. Prod. Res.}, 55(20):5902--5915, October 2017.

\bibitem[PCC18]{Pansart2018}
Lucie Pansart, Nicolas Catusse, and Hadrien Cambazard.
\newblock Exact algorithms for the order picking problem.
\newblock {\em Comput. Oper. Res.}, 100:117--127, December 2018.

\bibitem[PGBAAD24]{PGAD:OrderBatchingReview}
Eduardo~G. Pardo, Sergio Gil-Borrás, Antonio Alonso-Ayuso, and Abraham Duarte.
\newblock Order batching problems: Taxonomy and literature review.
\newblock {\em European Journal of Operational Research}, 313(1):1--24, 2024.

\bibitem[PN22]{Pinto2022-BatchingReview}
Anderson Rog{\'e}rio~Faia Pinto and Marcelo~Seido Nagano.
\newblock A comprehensive review of batching problems in low-level
  picker-to-parts systems with order due dates: Main gaps, trade-offs, and
  prospects for future research.
\newblock {\em J. Manuf. Syst.}, 65:1--18, October 2022.

\bibitem[RdK01]{Roodbergen2001}
Kees~Jan Roodbergen and Ren{\'e} de~Koster.
\newblock Routing order pickers in a warehouse with a middle aisle.
\newblock {\em Eur. J. Oper. Res.}, 133(1):32--43, August 2001.

\bibitem[RMSES25]{Rahmani_Mokarrari2025-ordersize}
Kimiya Rahmani~Mokarrari, Taraneh Sowlati, Jeffrey English, and Michael
  Starkey.
\newblock Optimization of warehouse picking to maximize the picked orders
  considering practical aspects.
\newblock {\em Appl. Math. Model.}, 137(115585):115585, January 2025.

\bibitem[RR83]{Ratliff1983}
H~Donald Ratliff and Arnon~S Rosenthal.
\newblock Order-picking in a rectangular warehouse: A solvable case of the
  traveling salesman problem.
\newblock {\em Oper. Res.}, 31(3):507--521, June 1983.

\bibitem[RV06]{Roodbergen2006-warehousedesign}
Kees~Jan Roodbergen and Iris F~A Vis.
\newblock A model for warehouse layout.
\newblock {\em IIE Trans.}, 38(10):799--811, October 2006.

\bibitem[SBK{\etalchar{+}}22]{Schiffer2022-ij}
Maximilian Schiffer, Nils Boysen, Patrick~S Klein, Gilbert Laporte, and Marco
  Pavone.
\newblock Optimal picking policies in e-commerce warehouses.
\newblock {\em Manage. Sci.}, 68(10):7497--7517, October 2022.

\bibitem[TADVN24]{Tran2024-routeassignment}
Son~Thai Tran, Rui~Jorge Almeida, Christof Defryn, and Inneke Van~Nieuwenhuyse.
\newblock Synchronous parallel heuristics for solving the joint order batching
  and picker routing problem.
\newblock In {\em 2024 {IEEE} Congress on Evolutionary Computation ({CEC})},
  pages 01--08. IEEE, June 2024.

\bibitem[TDK24]{Tutam2024-cobot}
Mahmut Tutam and Ren{\'e} De~Koster.
\newblock To walk or not to walk? designing intelligent order picking
  warehouses with collaborative robots.
\newblock {\em Transp. Res. Part E: Logist. Trans. Rev.}, 190(103696):103696,
  October 2024.

\bibitem[TWBT10]{Tompkins2010-FacilityPlanning}
James~A Tompkins, John~A White, Yavuz~A Bozer, and J~Tanchoco.
\newblock {\em Facilities Planning}.
\newblock Wiley-Blackwell, Hoboken, NJ, 4 edition, January 2010.

\bibitem[ZLF{\etalchar{+}}19]{Zhou2019-warehousedesign}
Li~Zhou, Jing Liu, Xiani Fan, Dongjie Zhu, Pingyu Wu, and Ning Cao.
\newblock Design of v-type warehouse layout and picking path model based on
  internet of things.
\newblock {\em IEEE Access}, 7:58419--58428, 2019.

\end{thebibliography}

\newpage
\begin{APPENDICES}

  \section{Additional Details}

  \subsection{Mathematical Model}\label{sec:app:mathmodel}
  In this section we give a full, self-contained statement of the mathematical model underlying the instances.

  Denote by $\Articles$ the set of articles and by $\Items$ the set of warehouse items.
  Each article $a\in \Articles$ has a volume $v_a \in \N$. Each $i\in \Items$ can fulfill its corresponding article $a_i
  \in \Articles$ and is located at a triple $(z_i, r_i, s_i)\in \Z^3$ corresponding to zone, rack, and aisle number.
  The set of orders is given by $\Orders$, where each order is given as a set of articles.

  Each picklist $p=(i_1, \dots, i_k) \in \Items^k$ consists of items within the same zone. For items in the same zone, we
  can define a distance function $d: \Items^2 \to \N$, which represents the length of the route in the warehouse model
  outlined above and we refer an interested reader to the implementation for an explicit definition. We assume the
  location of the depot in each zone $z\in \Zones$ to be at $(z, 0, 0)$ and alter the notation to define the distance from the
  depot as a function $d_0: \Items \to \N$ as the distance from the depot in the corresponding zone to the item. The length or \emph{cost} of a picklist, which is the main objective we optimize on, is thus given as
  \[ \ell(p) = d_0(i_1) + \sum_{j=1}^{k-1}d(i_j, i_{j+1}) + d_0(i_k). \]

  A collection of several picklists is called a \emph{batch} $b$.
  We do not require the items of a single picklist to correspond to a precise set of orders, as a picklist is restricted to a single zone and the corresponding items might have to be picked up across different zones.
  A batch, on the other hand, should always be comprised of picklists that pick up only those items that fulfill a certain set of orders, and no order should be fulfilled only partially.

  Additional parameters for the problem are the picker capacity $ \Capacity \in \N$, the item goal $ \Goal\in \N $, and a
  bound on the size of each batch $\Batch \in \N$.

  A feasible solution for the batching problem consists of a list of batches $b_1,\dots, b_k$ (where $k$ is not fixed) and
  fulfilled orders $O_1,\dots, O_k\subset \Orders$ with the following properties. Denote the set of feasible solutions as
  $\Feasible$. Remember that any sets are to be considered as multisets.
  \begin{enumerate}
    \item Every item is only picked at most once across all picklists:
      \[ \forall i\in \Items: \card{\{j \mid j = i, j\in p, p\in (b_1\cup \dots \cup b_k) \}} \leq 1 \]
    \item The total item goal is fulfilled:
      \[ \sum_{j=1}^k \sum_{p\in b_j} \card{p} \geq \Goal. \]
    \item For each $j \in 1,\dots, k$ we require:
      \begin{enumerate}
        \item The total volume of items in each picklist in the batch does not exceed $\Capacity$:
          \[ \forall p\in b_j: \sum_{i\in p}v_{a_i} \leq \Capacity. \]
        \item The batch size is respected: $\card{O_j} \leq \Batch$.
        \item The requested orders for each batch are precisely fulfilled, i.e.\ the number of times each article is ordered corresponds to the articles of all picklists in the batch:
          \[ \bigcup_{o\in O_j} o = \{a_i \mid i\in p, p\in b_j \}. \]
      \end{enumerate}
  \end{enumerate}

  The objective function to be minimized is the total length of picklists among all feasible solutions:
  \[ \min_{(b_1,\dots, b_k)\in \Feasible} \sum_{j=1}^{k} \sum_{p\in b_j} \ell(p). \]

  \subsection{FDGA with Strict Age Objective}\label{sec:app:fdgastrictage}

  While the previous approach accomplishes the completion of old orders implicitly due to the modification of the costs, we have also devised a second approach that gives explicit hard guarantees by design of the algorithm.

  This variant of \fdga takes as input a list of non-increasing order numbers $B = b_0,\dots, b_m$, where $b_m = 0$.
  We also call the $b_i$ \emph{bucket size of age $i$} or \emph{completion targets}.
  Each $b_i$ indicates how many orders of age $i$ may be carried over to the next phase of the algorithm.
  Thus, no order can exceed age $m$ over the course of the algorithm.

  As an example, consider $B = (1000,\ 500,\allowbreak 250,\ 0)$.
  The algorithm can carry over 1000 new orders, 500 orders of age 1 and 250 orders of age 2 to the next phase.
  Any order older than age 2 \emph{must} be completed in the given phase.

  Clearly, these completion targets may interfere with the item goal if not appropriately set.
  When the number of orders (of any given age) in a phase vastly exceeds the bucket size, all but some of the orders have to be completed in this phase and the corresponding number of items may have to exceed the item goal substantially.
  For a regular warehouse setting with an approximately constant number $n$ of incoming orders in each phase, it suffices to have $b_0 \geq n$.
  The number of orders of age $i\geq 1$ to be fulfilled in any given phase is no larger than $ b_{i-1}-b_{i}$, and therefore the total number is no larger than $n-b_0 + \sum_{i=1}^{m} b_{i-1}-b_i = n$.
  Thus, if the warehouse manages to handle the incoming order stream, we can also apply this procedure without having to exceed the item goal in any way.

  As a remark for the experimental implementation with synthetic data, we note that we start the initial phase with a certain amount of pre-existing orders in the warehouse.
  These orders are given artificial ages in the initial phase according to a backfill procedure which assigns ages to them such that the age buckets are entirely filled from the back. In the example above for 600 initial orders, we would have 250 randomly selected orders of age 2 and 350 of age 1.
  In order for this backfill procedure to be applied, the sum of bucket sizes must exceed the number of initial orders.
  After at most $m$ phases, all of the initial orders are completed and the order ages converge to the ``natural'' distribution of the algorithm.

  Given completion targets $B$, the algorithm computes the order selection in the following way:

  \begin{enumerate}
    \item From the average order size, estimate the number of batches $b$ required to approximately fulfill the completion target.
    \item From each set of orders of a given age, randomly sample as many orders as necessary to achieve the completion targets.
      These orders form the set $\mathcal{M}$ of \emph{mandatory orders}.
    \item Within each batch, we use the \fdga procedure on $\mathcal{M}$ to find the first $|\mathcal{M}|/b$ orders.
      Thereby, we attempt to evenly spread the mandatory orders across all batches.
    \item Then, we use \fdga on the full order set to complete the batch up to the maximum batch size, i.e.\ we use \emph{both} $\mathcal{M}$ and the remaining orders that are not selected as mandatory.
    \item If no more mandatory orders are left during the first phase, we simply use the remaining orders to compute a standard batch with \fdga.
  \end{enumerate}

  In our experiments, we use the following implementations of the bucket sizes. Note that both of them implicitly use (an estimate of) the number $n$ of newly added orders in each phase.

  \begin{itemize}
    \item \fdgalinearn{$m$} implements a linear decline of bucket sizes from $n$ to 0 at age $m$.
      The buckets are thus given by $b_i = \max(0, \lceil n\cdot (m-i)/m\rceil)$
    \item \fdgacutoffn{$m$} does not impose any restrictions on the orders up to age $m$, but then requires all orders of that age to be completed.
      Thus we have $b_0,\dots, b_{m-1} = n$ and $b_m = 0$.
  \end{itemize}

  \section{Instance Parameters}\label{sec:app:instanceparams}

  Before we present our data, we outline some details of the computation that are not addressed in the abstract problem model in \cref{sec:singlephase-model}.
  Our computational experiments and reference algorithms (\dga and \rdga as explained in the paper) are based on \cite{ZalandoPaper} and their public code repository for this project (\url{https://github.com/zalandoresearch/batching-benchmarks/}).
  To allow for better comparability with their results, we maintain their instance generation parameters for the different categories, which we reproduce in \cref{tab:instanceparams}.

  The complete generation procedure of \cite{ZalandoPaper} is not described in their article, but can only be found in their code repository, and is summarized in the following:
  \begin{itemize}
    \item The number of different articles in the warehouse is approximately $\frac13$ of the number of warehouse items.
      Each of these articles is assigned a random article volume taken as $\max(1, \operatorname{round}(V))$, where $V\sim
      \Gamma(2, 20)$ is drawn from a gamma distribution and $\operatorname{round}(x)$ denotes rounding $x$ to the nearest integer value.
    \item Then, approximately three items are generated for each article. The remainder of the
      warehouse is filled with items that correspond to uniformly sampled articles. The location of each item is uniformly sampled across the available zones, rows and aisles.
    \item The orders are finally generated from the list of available items. The order size is first sampled from a
      weighted distribution of sizes between 2 and 6 and then the articles of randomly chosen items are assigned. In order
      to ensure that every order can be completed, the chosen items are removed from the list for the generation of the
      next order.
  \end{itemize}

  \newpage
  \section{Additional Data: Singlephase Instances in Different Variants}\label{sec:app:singlephasedata}

  \subsection{medium}

  \begin{table}[H] 
    \TABLE{Comparing the heuristic routing subroutine with the improved CVRP-based solution for medium instances.
    \label{tab:singlephasecvrp-medium}}
    {
      \begin{tabular}{c@{\hskip 2\tabcolsep}|cc|ccc}
        \toprule
        & \multicolumn{2}{c|}{greedy routing} & \multicolumn{2}{c}{CVRP} &  \\
        Algorithm  & Cost    & Time   & Cost             & Time & Relative Performance \\
        \midrule
        \fdga & 62354.4 & 8.7s & 61526.4 & 9.7s & -1.3\% \\
        \fdgaskip & 62203.2 & 7.9s & 61215.6 & 8.9s & -1.6\% \\
        \ldga & 70684.4 & 5.7s & 68294.4 & 6.6s & -3.4\% \\
        \rdgasqrt & 76453.2 & 1.8s & 75569.2 & 2.7s & -1.2\% \\
        \rdgan{50} & 77449.2 & 1.4s & 77311.6 & 2.3s & -0.2\% \\
        \rdga & 139287.2 & 0.2s & 135510.8 & 1.0s & -2.7\% \\
        \bottomrule
    \end{tabular}}{All values are averaged over five different instances and rounded to a tenth.}
  \end{table}

  \subsection{medium-fewzones}

  \begin{table}[H]
    \TABLE{Further computational measures for medium instances with few and fewer zones.
    \label{tab:singlephaseinsights-fewerzones}}
    {\footnotesize
      \begin{tabular}{c@{\hskip 2\tabcolsep}|ccc|ccc}\toprule
        & \multicolumn{3}{c|}{medium-fewzones} & \multicolumn{3}{c}{medium-fewerzones} \\
        Algorithm & \shortstack{Item distance\\from conveyor}  & Order size & Picklist length & \shortstack{Item distance\\from conveyor}  & Order size & Picklist length \\\midrule
        \fdga & 24.8 & 2.2 &  10.3 & 26.1 & 2.2  & 18.6 \\
        \fdgaskip & 24.9 & 2.2 &  10.3 & 26.1 & 2.2  & 18.6 \\
        \ldga & 25.0 & 2.3 &  10.8& 26.5 & 2.3  & 18.7 \\
        \rdgan{50} & 23.7 & 2.4  & 11.2 & 26.0 & 2.4  & 18.9\\
        \rdgasqrt & 23.6 & 2.3 & 10.9 & 25.7 & 2.3 & 19.0\\
        \rdga & 36.1 & 2.7  & 13.0& 38.9 & 2.6  & 19.2 \\
        (Instance) & 33.93 & 2.65 & $-$ &  34.01 & 2.64  &$-$\\\bottomrule
    \end{tabular}}{All values are averaged over five different instances and represent the mean value over all batched items, picklists, and orders respectively. The reference values refer to all orders/items in the instance as compared to the ones picked by the algorithm. The picklist lengths correspond to average picker loads of 40\% and 70-80\%, respectively.}
  \end{table}

  \begin{table}[H]
    \TABLE{Comparing the heuristic routing subroutine with the improved CVRP-based solution for medium-fewzones instances.
    \label{tab:singlephasecvrp-medium-fewzones}}
    {
      \begin{tabular}{c@{\hskip 2\tabcolsep}|cc|ccc}
        \toprule
        & \multicolumn{2}{c|}{greedy routing} & \multicolumn{2}{c}{CVRP} &  \\
        Algorithm  & Cost    & Time   & Cost             & Time & Relative Performance \\
        \midrule
        \dga & 41446.4 & 185.5s & 39017.2 & 191.2s & -5.9\% \\
        \fdga & 41446.4 & 30.4s & 39017.2 & 34.4s & -5.9\% \\
        \fdgaskip & 41682.4 & 27.7s & 38881.2 & 29.8s & -6.7\% \\
        \ldga & 43815.6 & 9.1s & 41380.0 & 10.7s & -5.6\% \\
        \ldga-global & 44480.4 & 17.7s & 41004.8 & 19.3s & -7.8\% \\
        \rdga & 102917.6 & 0.2s & 91208.8 & 1.5s & -11.4\% \\
        \rdgan{50} & 53383.2 & 2.3s & 51147.2 & 3.7s & -4.2\% \\
        \rdgasqrt & 51320.4 & 2.9s & 49572.8 & 4.4s & -3.4\% \\
        \bottomrule
    \end{tabular}}{All values are averaged over five different instances and rounded to a tenth.}
  \end{table}

  \subsection{medium-fewerzones}

  \begin{table}[H]
    \TABLE{Comparing the heuristic routing subroutine with the improved CVRP-based solution for medium-fewerzones instances.
    \label{tab:singlephasecvrp-medium-fewerzones}}
    {
      \begin{tabular}{c@{\hskip 2\tabcolsep}|cc|ccc}
        \toprule
        & \multicolumn{2}{c|}{greedy routing} & \multicolumn{2}{c}{CVRP} &  \\
        Algorithm  & Cost    & Time   & Cost             & Time & Relative Performance \\
        \midrule
        \fdga & 29281.6 & 74.7s & 26103.2 & 82.8s & -10.9\% \\
        \fdgaskip & 28498.4 & 67.7s & 25474.0 & 72.9s & -10.6\% \\
        \ldga & 33404.4 & 21.1s & 29532.8 & 30.4s & -11.6\% \\
        \rdgasqrt & 36911.6 & 6.1s & 34086.0 & 11.6s & -7.7\% \\
        \rdgan{50} & 38494.0 & 4.7s & 35381.2 & 10.0s & -8.1\% \\
        \rdga & 83658.4 & 0.3s & 64500.4 & 6.2s & -22.9\% \\
        \bottomrule
    \end{tabular}}{All values are averaged over five different instances and rounded to a tenth.}
  \end{table}

  \subsection{medium-longbatch}

  \begin{table}[H]
    \TABLE{Benchmark data of selected algorithms for medium-longbatch instances. \label{tab:singlephaseperformance-medium-longbatch}}
    {
      \begin{tabular}{c@{\hskip 2\tabcolsep}|ccr}\toprule
        Algorithm  & Cost & Relative Performance & Time  \\\midrule
        \dga & 57119.6 & +0.3\% & 200.6s \\
        \fdga & 57119.6 & +0.3\% & 7.6s \\
        \fdgaskip & 56940.0 & (reference) & 7.1s \\
        \ldga & 65272.8 & +14.6\% & 9.1s \\
        \ldga-global & 59092.8 & +3.8\% & 9.8s \\
        \rdga & 104873.2 & +84.2\% & 0.3s \\
        \rdgan{50} & 63814.0 & +12.1\% & 2.6s \\
        \rdgasqrt & 62356.0 & +9.5\% & 3.3s \\
        \bottomrule
    \end{tabular}}{All values are averaged over five different medium-longbatch instances and rounded to the nearest tenth.}
  \end{table}

  \begin{table}[H]
    \TABLE{Further computational measures for medium-longbatch instances.
    \label{tab:singlephaseinsights-medium-longbatch}}
    {
      \begin{tabular}{c@{\hskip 2\tabcolsep}|ccc}\toprule
        Algorithm & \shortstack{Item distance\\from conveyor}  & Order size & Picklist length \\\midrule
        \dga & 26.2 & 2.4 & 10.8 \\
        \fdga & 26.2 & 2.4 &  10.8 \\
        \fdgaskip & 26.1 & 2.4 &  10.8 \\
        \ldga & 27.1 & 2.5 &  11.6 \\
        \ldga-global & 26.0 & 2.4 & 11.0 \\
        \rdga & 35.8 & 2.7 &  12.9 \\
        \rdgan{50} & 24.9 & 2.4 & 11.1 \\
        \rdgasqrt & 24.7 & 2.4 &  11.1 \\
        (Instance) & 34.03 & 2.65 & $-$ \\
        \bottomrule
    \end{tabular}}{All values are averaged over five different medium-longbatch instances and represent the mean value over all batched items, picklists, and orders respectively. The reference values refer to all orders/items in the instance as compared to the ones picked by the algorithm.}
  \end{table}

  \begin{table}[H]
    \TABLE{Comparing the heuristic routing subroutine with the improved CVRP-based solution for medium-longbatch instances.
    \label{tab:singlephasecvrp-medium-longbatch}}
    {
      \begin{tabular}{c@{\hskip 2\tabcolsep}|cc|ccc}
        \toprule
        & \multicolumn{2}{c|}{greedy routing} & \multicolumn{2}{c}{CVRP} &  \\
        Algorithm  & Cost    & Time   & Cost             & Time & Relative Performance \\
        \midrule
        \dga & 57119.6 & 200.6s & 53829.6 & 202.7s & -5.8\% \\
        \fdga & 57119.6 & 7.6s & 53829.6 & 9.6s & -5.8\% \\
        \fdgaskip & 56940.0 & 7.1s & 53597.6 & 8.8s & -5.9\% \\
        \ldga & 65272.8 & 9.1s & 60786.8 & 10.6s & -6.9\% \\
        \ldga-global & 59092.8 & 9.8s & 56408.8 & 11.7s & -4.5\% \\
        \rdga & 104873.2 & 0.3s & 102068.8 & 1.7s & -2.7\% \\
        \rdgan{50} & 63814.0 & 2.6s & 61192.0 & 4.3s & -4.1\% \\
        \rdgasqrt & 62356.0 & 3.3s & 59592.8 & 4.9s & -4.4\% \\
        \bottomrule
    \end{tabular}}{All values are averaged over five different instances and rounded to a tenth.}
  \end{table}

  \subsection{large}

  \begin{table}[H]
    \TABLE{Comparing the heuristic routing subroutine with the improved CVRP-based solution for large instances.
    \label{tab:singlephasecvrp-large}}
    {
      \begin{tabular}{c@{\hskip 2\tabcolsep}|cc|ccc}
        \toprule
        & \multicolumn{2}{c|}{greedy routing} & \multicolumn{2}{c}{CVRP} &  \\
        Algorithm  & Cost    & Time   & Cost             & Time & Relative Performance \\
        \midrule
        \dga & 524925.2 & 10548.6s & 515961.6 & 10556.8s & -1.7\% \\
        \fdga & 524925.2 & 591.8s & 515961.6 & 603.6s & -1.7\% \\
        \fdgaskip & 526366.8 & 545.5s & 517194.4 & 557.3s & -1.7\% \\
        \ldga & 597082.8 & 113.8s & 590370.4 & 123.4s & -1.1\% \\
        \ldga-global & 583514.4 & 198.4s & 574965.6 & 207.1s & -1.5\% \\
        \rdga & 1489700.8 & 9.8s & 1479903.2 & 19.7s & -0.7\% \\
        \rdgan{50} & 843182.0 & 22.5s & 837952.4 & 31.5s & -0.6\% \\
        \rdgasqrt & 734645.2 & 59.1s & 730166.0 & 68.7s & -0.6\% \\
        \bottomrule
    \end{tabular}}{All values are averaged over five different instances and rounded to a tenth.}
  \end{table}

  \subsection{large-fewzones}

  \begin{table}[H]
    \TABLE{Benchmark data of selected algorithms for large-fewzones instances. \label{tab:singlephaseperformance-large-fewzones}}
    {
      \begin{tabular}{c@{\hskip 2\tabcolsep}|ccr}\toprule
        Algorithm  & Cost & Relative Performance & Time  \\\midrule
        \dga & 293734.0 & +0.2\% & 19024.2s \\
        \fdga & 293734.0 & +0.2\% & 3142.6s \\
        \fdgaskip & 293233.6 & (reference) & 2841.0s \\
        \ldga & 402608.4 & +37.3\% & 569.5s \\
        \ldga-global & 388380.8 & +32.4\% & 1371.5s \\
        \rdga & 1003955.2 & +242.4\% & 9.4s \\
        \rdgan{50} & 513872.4 & +75.2\% & 31.1s \\
        \rdgasqrt & 432059.6 & +47.3\% & 100.6s \\
        \bottomrule
    \end{tabular}}{All values are averaged over five different large-fewzones instances and rounded to the nearest tenth.}
  \end{table}

  \begin{table}[H]
    \TABLE{Further computational measures for large-fewzones instances.
    \label{tab:singlephaseinsights-large-fewzones}}
    {
      \begin{tabular}{c@{\hskip 2\tabcolsep}|ccc}\toprule
        Algorithm & \shortstack{Item distance\\from conveyor}  & Order size & Picklist length \\\midrule
        \dga & 24.8 & 2.2 & 10.0 \\
        \fdga & 24.8 & 2.2  & 10.0 \\
        \fdgaskip & 24.7 & 2.2 &  10.0 \\
        \ldga & 25.5 & 2.3 &  10.5 \\
        \ldga-global & 24.7 & 2.3 &10.5 \\
        \rdga & 35.8 & 2.6 &  12.8 \\
        \rdgan{50} & 23.7 & 2.4 &  11.2 \\
        \rdgasqrt & 23.7 & 2.3 & 10.6 \\
        (Instance) & 34.02 & 2.64 & $-$ \\
        \bottomrule
    \end{tabular}}{All values are averaged over five different large-fewzones instances and represent the mean value over all batched items, picklists, and orders respectively. The reference values refer to all orders/items in the instance as compared to the ones picked by the algorithm.}
  \end{table}

  \begin{table}[H]
    \TABLE{Comparing the heuristic routing subroutine with the improved CVRP-based solution for large-fewzones instances.
    \label{tab:singlephasecvrp-large-fewzones}}
    {
      \begin{tabular}{c@{\hskip 2\tabcolsep}|cc|ccc}
        \toprule
        & \multicolumn{2}{c|}{greedy routing} & \multicolumn{2}{c}{CVRP} &  \\
        Algorithm  & Cost    & Time   & Cost             & Time & Relative Performance \\
        \midrule
        \dga & 293734.0 & 19024.2s & 273636.0 & 18974.8s & -6.8\% \\
        \fdga & 293734.0 & 3142.6s & 273649.6 & 3136.7s & -6.8\% \\
        \fdgaskip & 293233.6 & 2841.0s & 274330.8 & 2861.8s & -6.4\% \\
        \ldga & 402608.4 & 569.5s & 374332.4 & 587.8s & -7.0\% \\
        \ldga-global & 388380.8 & 1371.5s & 366252.0 & 1384.7s & -5.7\% \\
        \rdga & 1003955.2 & 9.4s & 884936.4 & 21.7s & -11.9\% \\
        \rdgan{50} & 513872.4 & 31.1s & 492736.0 & 45.4s & -4.1\% \\
        \rdgasqrt & 432059.6 & 100.6s & 413892.8 & 119.9s & -4.2\% \\
        \bottomrule
    \end{tabular}}{All values are averaged over five different instances and rounded to a tenth.}
  \end{table}

  \subsection{large-fewerzones}

  \begin{table}[H]
    \TABLE{Benchmark data of selected algorithms for large-fewerzones instances. \label{tab:singlephaseperformance-large-fewerzones}}
    {
      \begin{tabular}{c@{\hskip 2\tabcolsep}|ccr}\toprule
        Algorithm  & Cost & Relative Performance & Time  \\\midrule
        \dga & 198821.6 & (reference) & 39523.7s \\
        \fdga & 198824.8 & +0.0\% & 7349.1s \\
        \fdgaskip & 199959.2 & +0.6\% & 7387.5s \\
        \ldga & 311075.2 & +56.5\% & 1211.4s \\
        \rdgasqrt & 298406.0 & +50.1\% & 199.0s \\
        \rdgan{50} & 369257.2 & +85.7\% & 55.2s \\
        \rdga & 826112.8 & +315.5\% & 9.8s \\
        \bottomrule
    \end{tabular}}{All values are averaged over five different large-fewerzones instances and rounded to the nearest tenth.}
  \end{table}

  \begin{table}[H]
    \TABLE{Further computational measures for large-fewerzones instances.
    \label{tab:singlephaseinsights-large-fewerzones}}
    {
      \begin{tabular}{c@{\hskip 2\tabcolsep}|ccc}\toprule
        Algorithm & \shortstack{Item distance\\from conveyor}  & Order size & Picklist length \\\midrule
        \dga & 25.7 & 2.1 & 18.0 \\
        \fdga & 25.7 & 2.1  & 18.0 \\
        \fdgaskip & 25.6 & 2.1  & 18.1 \\
        \ldga & 26.5 & 2.3  & 18.7 \\
        \rdgasqrt & 24.8 & 2.3  & 18.7 \\
        \rdgan{50} & 25.3 & 2.3& 18.8 \\
        \rdga & 38.4 & 2.6 &  19.3 \\
        (Instance) & 34.01 & 2.64 & $-$ \\
        \bottomrule
    \end{tabular}}{All values are averaged over five different large-fewerzones instances and represent the mean value over all batched items, picklists, and orders respectively. The reference values refer to all orders/items in the instance as compared to the ones picked by the algorithm.}
  \end{table}

  \begin{table}[H]
    \TABLE{Comparing the heuristic routing subroutine with the improved CVRP-based solution for large-fewerzones instances.
    \label{tab:singlephasecvrp-large-fewerzones}}
    {
      \begin{tabular}{c@{\hskip 2\tabcolsep}|cc|ccc}
        \toprule
        & \multicolumn{2}{c|}{greedy routing} & \multicolumn{2}{c}{CVRP} &  \\
        Algorithm  & Cost    & Time   & Cost             & Time & Relative Performance \\
        \midrule
        \dga & 198821.6 & 39523.7s & 175403.6 & 39742.8s & -11.8\% \\
        \fdga & 198824.8 & 7349.1s & 175394.4 & 7408.7s & -11.8\% \\
        \fdgaskip & 199959.2 & 7387.5s & 176067.6 & 7628.1s & -11.9\% \\
        \ldga & 311075.2 & 1211.4s & 273192.0 & 1277.7s & -12.2\% \\
        \rdgasqrt & 298406.0 & 199.0s & 272717.6 & 263.0s & -8.6\% \\
        \rdgan{50} & 369257.2 & 55.2s & 337725.6 & 109.6s & -8.5\% \\
        \rdga & 826112.8 & 9.8s & 640261.2 & 81.5s & -22.5\% \\
        \bottomrule
    \end{tabular}}{All values are averaged over five different instances and rounded to a tenth.}
  \end{table}

  \subsection{large-longbatch}

  \begin{table}[H]
    \TABLE{Benchmark data of selected algorithms for large-longbatch instances. \label{tab:singlephaseperformance-large-longbatch}}
    {
      \begin{tabular}{c@{\hskip 2\tabcolsep}|ccr}\toprule
        Algorithm  & Cost & Relative Performance & Time  \\\midrule
        \dga & 441973.2 & +0.1\% & 14517.6s \\
        \fdga & 441973.2 & +0.1\% & 402.4s \\
        \fdgaskip & 441499.6 & (reference) & 375.0s \\
        \ldga & 498312.4 & +12.9\% & 127.2s \\
        \rdgasqrt & 534064.4 & +21.0\% & 80.2s \\
        \rdgan{50} & 624935.6 & +41.5\% & 27.6s \\
        \rdga & 1164882.0 & +163.8\% & 10.0s \\
        \bottomrule
    \end{tabular}}{All values are averaged over five different large-longbatch instances and rounded to the nearest tenth.}
  \end{table}

  \begin{table}[H]
    \TABLE{Further computational measures for large-longbatch instances.
    \label{tab:singlephaseinsights-large-longbatch}}
    {
      \begin{tabular}{c@{\hskip 2\tabcolsep}|ccc}\toprule
        Algorithm & \shortstack{Item distance\\from conveyor}  & Order size & Picklist length \\\midrule
        \dga & 24.4 & 2.2  & 6.1 \\
        \fdga & 24.4 & 2.2  & 6.1 \\
        \fdgaskip & 24.4 & 2.2 & 6.1 \\
        \ldga & 24.5 & 2.3  & 6.2 \\
        \rdgasqrt & 23.1 & 2.3  & 6.1 \\
        \rdgan{50} & 22.9 & 2.4 & 6.2 \\
        \rdga & 34.7 & 2.6  & 6.6 \\
        (Instance) & 34.03 & 2.64 & $-$ \\
        \bottomrule
    \end{tabular}}{All values are averaged over five different large-longbatch instances and represent the mean value over all batched items, picklists, and orders respectively. The reference values refer to all orders/items in the instance as compared to the ones picked by the algorithm.}
  \end{table}

  \begin{table}[H]
    \TABLE{Comparing the heuristic routing subroutine with the improved CVRP-based solution for large-longbatch instances.
    \label{tab:singlephasecvrp-large-longbatch}}
    {
      \begin{tabular}{c@{\hskip 2\tabcolsep}|cc|ccc}
        \toprule
        & \multicolumn{2}{c|}{greedy routing} & \multicolumn{2}{c}{CVRP} &  \\
        Algorithm  & Cost    & Time   & Cost             & Time & Relative Performance \\        \midrule
        \dga & 441973.2 & 14517.6s & 425062.8 & 14444.7s & -3.8\% \\
        \fdga & 441973.2 & 402.4s & 425062.8 & 418.5s & -3.8\% \\
        \fdgaskip & 441499.6 & 375.0s & 425532.8 & 390.8s & -3.6\% \\
        \ldga & 498312.4 & 127.2s & 478226.4 & 137.2s & -4.0\% \\
        \rdgasqrt & 534064.4 & 80.2s & 520688.0 & 91.5s & -2.5\% \\
        \rdgan{50} & 624935.6 & 27.6s & 613729.2 & 36.8s & -1.8\% \\
        \rdga & 1164882.0 & 10.0s & 1079676.8 & 17.9s & -7.3\% \\
        \bottomrule
    \end{tabular}}{All values are averaged over five different instances and rounded to a tenth.}
  \end{table}

  \section{Additional Data: Multiphase Instances in Different Variants}\label{sec:app:multiphasedata}

  In this section, we augment the discussion from the main text by providing additional numerical results for the medium and large instance variants not outlined there (and including the data already provided in the main text).

  \begin{table}[H] 
    \TABLE{Benchmark data comparing different multi-phase algorithms for medium instances. \label{tab:app:multiphase-performance-medium}}
    {
      \begin{tabular}{c@{\hskip 2\tabcolsep}|cccc}\toprule
        Algorithm & Total cost & Rel.\ cost & \shortstack{Avg.\ cost\\ per item} & \shortstack{Avg.\ cost\\ per order} \\\midrule
        \fdga & 9120363 & (reference) & 34.94 & 90.98 \\
        \fdgacutoffn{5} & 11291446 & $+23.8\% $ & 43.85 & 115.47 \\
        \fdgacutoffn{10} & 10673561 & $+17.0\% $ & 41.32 & 108.36 \\
        \fdgacutoffn{15} & 10173553 & $+11.5\% $ & 39.34 & 102.96 \\
        \fdgacutoffn{20} & 9893970 & $+8.5\% $ & 38.19 & 99.85 \\
        \fdgacutoffn{25} & 9728927 & $+6.7\% $ & 37.48 & 97.89 \\
        \fdgacutoffn{30} & 9600925 & $+5.3\% $ & 36.99 & 96.58 \\
        \fdgalinearn{10} & 10900664 & $+19.5\% $ & 42.07 & 110.38 \\
        \fdgalinearn{15} & 10290667 & $+12.8\% $ & 39.71 & 103.95 \\
        \fdgalinearn{20} & 9941275 & $+9.0\% $ & 38.36 & 100.31 \\
        \fdgalinearn{25} & 9752540 & $+6.9\% $ & 37.6 & 98.21 \\
        \fdgalinearn{30} & 9650298 & $+5.8\% $ & 37.09 & 96.82 \\
        \fdgasoftg{0.8} & 9307641 & $+2.1\% $ & 35.95 & 94.61 \\
        \fdgasoftg{0.9} & 8840642 & $-3.1\% $ & 34.16 & 89.73 \\
        \fdgasoftg{0.95} & 8697709 & $-4.6\% $ & 33.55 & 87.98 \\
        \fdgasoftg{0.99} & 8849031 & $-3.0\% $ & 34.02 & 88.89 \\
        \ldga & 9243751 & $+1.4\% $ & 35.61 & 93.16 \\
        \rdgasqrt & 10826678 & $+18.7\% $ & 41.6 & 108.8 \\\bottomrule
    \end{tabular}}{The parameter for the number of new orders is set to 1000 if applicable and approximately corresponds to the number of new orders in each phase. The column for relative cost compares to the reference of \fdga. The values are averaged over 5 instances and rounded to appropriate precision.  }
  \end{table}

  \begin{table}[H] 
    \TABLE{Further computational measures for different algorithms on medium instances.

    \label{tab:multiphaseinsights-medium-new}}
    {
      \begin{tabular}{c@{\hskip 2\tabcolsep}|cccc}\toprule
        Algorithm & \shortstack{Average\\ order length} & \shortstack{Maximum\\ order age} & \shortstack{Average\\ order age} & \shortstack{Median\\ order age} \\\midrule
        \fdga & 3.52 & 100.0 & 16.88 & 4.8 \\
        \fdgacutoffn{5} & 2.71 & 5.0 & 2.34 & 2.0 \\
        \fdgacutoffn{10} & 2.96 & 10.0 & 3.63 & 3.0 \\
        \fdgacutoffn{15} & 3.1 & 15.0 & 4.67 & 3.0 \\
        \fdgacutoffn{20} & 3.18 & 20.0 & 5.58 & 3.2 \\
        \fdgacutoffn{25} & 3.25 & 25.0 & 6.54 & 4.0 \\
        \fdgacutoffn{30} & 3.28 & 30.0 & 7.25 & 4.0 \\
        \fdgalinearn{10} & 2.93 & 10.0 & 3.42 & 3.0 \\
        \fdgalinearn{15} & 3.08 & 15.0 & 4.58 & 3.0 \\
        \fdgalinearn{20} & 3.16 & 20.0 & 5.52 & 3.2 \\
        \fdgalinearn{25} & 3.23 & 25.0 & 6.43 & 4.0 \\
        \fdgalinearn{30} & 3.29 & 30.0 & 7.23 & 4.0 \\
        \fdgasoftg{0.8} & 2.74 & 7.6 & 2.22 & 2.0 \\
        \fdgasoftg{0.9} & 2.85 & 10.8 & 2.53 & 2.0 \\
        \fdgasoftg{0.95} & 2.95 & 17.2 & 3.17 & 2.0 \\
        \fdgasoftg{0.99} & 3.23 & 48.8 & 6.0 & 3.0 \\
        \ldga & 3.12 & 100.0 & 11.52 & 3.2 \\
        \rdgasqrt & 3.15 & 100.0 & 13.66 & 4.0 \\\bottomrule
    \end{tabular}}{All values are averaged over five different instances.
    The values are computed for the set of orders in the \emph{last} phase of each run.}
  \end{table}

  \begin{figure}[H]
    \FIGURE{
    \includegraphics[width=0.8\textwidth]{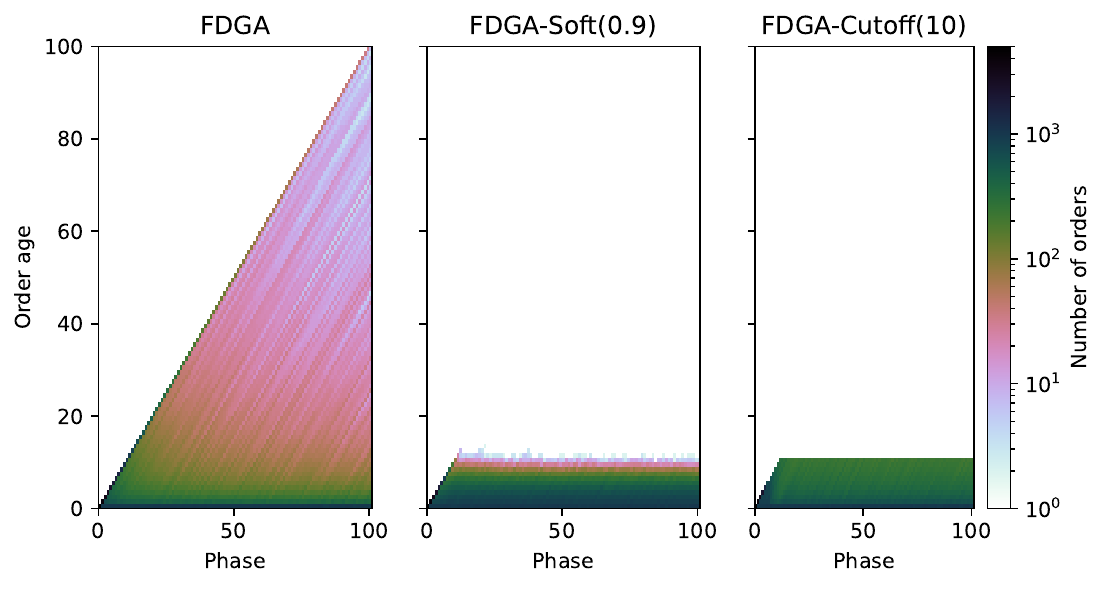}}
    {Age distribution data of single runs of \fdga, \fdgasoftg{0.9}, and \fdgacutoffn{10}, respectively, on a medium instance.
      While many orders are completed right away, many of them remain in the problem until the end of the experiment.
    \label{fig:app:fdga-age}}{}
  \end{figure}

  \begin{table}[H] 
    \TABLE{Benchmark data comparing different multi-phase algorithms for medium-fewzones instances.   \label{tab:app:multiphase-performance-medium-fewzones}}
    {
      \begin{tabular}{c@{\hskip 2\tabcolsep}|cccc}\toprule      Algorithm & Total cost & Rel.\ cost & \shortstack{Avg.\ cost\\ per item} & \shortstack{Avg.\ cost\\ per order} \\\midrule
        \fdga & 5609744 & (reference) & 21.51 & 55.97 \\
        \fdgacutoffn{5} & 7023362 & $+25.2\% $ & 27.28 & 71.84 \\
        \fdgacutoffn{10} & 6652562 & $+18.6\% $ & 25.75 & 67.46 \\
        \fdgacutoffn{15} & 6352230 & $+13.2\% $ & 24.59 & 64.29 \\
        \fdgacutoffn{20} & 6188696 & $+10.3\% $ & 23.87 & 62.34 \\
        \fdgacutoffn{25} & 6067686 & $+8.2\% $ & 23.4 & 61.07 \\
        \fdgacutoffn{30} & 5991385 & $+6.8\% $ & 23.04 & 60.12 \\
        \fdgalinearn{10} & 6750239 & $+20.3\% $ & 26.15 & 68.55 \\
        \fdgalinearn{15} & 6423991 & $+14.5\% $ & 24.81 & 64.88 \\
        \fdgalinearn{20} & 6226611 & $+11.0\% $ & 24.03 & 62.77 \\
        \fdgalinearn{25} & 6113034 & $+9.0\% $ & 23.53 & 61.43 \\
        \fdgalinearn{30} & 6018649 & $+7.3\% $ & 23.16 & 60.43 \\
        \fdgasoftg{0.8} & 5539346 & $-1.3\% $ & 21.39 & 56.25 \\
        \fdgasoftg{0.9} & 5315791 & $-5.2\% $ & 20.55 & 53.93 \\
        \fdgasoftg{0.95} & 5299643 & $-5.5\% $ & 20.43 & 53.51 \\
        \fdgasoftg{0.99} & 5469906 & $-2.5\% $ & 21.01 & 54.84 \\
        \ldga & 5566124 & $-0.8\% $ & 21.41 & 55.9 \\
        \rdgasqrt & 6620678 & $+18.0\% $ & 25.44 & 66.5 \\\bottomrule
    \end{tabular}}{The parameter for the number of new orders is set to 1000 if applicable and approximately corresponds to the number of new orders in each phase. The column for relative cost compares to the reference of \fdga. The values are averaged over 5 instances and rounded to appropriate precision.}
  \end{table}

  \begin{table}[H] 
    \TABLE{Further computational measures for different algorithms on medium-fewzones instances.

    \label{tab:multiphaseinsights-medium-fewzones}}
    {
      \begin{tabular}{c@{\hskip 2\tabcolsep}|cccc}\toprule
        Algorithm & \shortstack{Average\\ order length} & \shortstack{Maximum\\ order age} & \shortstack{Average\\ order age} & \shortstack{Median\\ order age} \\\midrule
        \fdga & 3.6 & 100.0 & 15.72 & 4.2 \\
        \fdgacutoffn{5} & 2.73 & 5.0 & 2.34 & 2.0 \\
        \fdgacutoffn{10} & 3.03 & 10.0 & 3.53 & 3.0 \\
        \fdgacutoffn{15} & 3.17 & 15.0 & 4.56 & 3.0 \\
        \fdgacutoffn{20} & 3.26 & 20.0 & 5.46 & 3.0 \\
        \fdgacutoffn{25} & 3.32 & 25.0 & 6.27 & 3.2 \\
        \fdgacutoffn{30} & 3.38 & 30.0 & 7.05 & 3.4 \\
        \fdgalinearn{10} & 3.0 & 10.0 & 3.37 & 3.0 \\
        \fdgalinearn{15} & 3.16 & 15.0 & 4.47 & 3.0 \\
        \fdgalinearn{20} & 3.26 & 20.0 & 5.43 & 3.0 \\
        \fdgalinearn{25} & 3.33 & 25.0 & 6.31 & 3.2 \\
        \fdgalinearn{30} & 3.38 & 30.0 & 7.11 & 3.8 \\
        \fdgasoftg{0.8} & 2.81 & 8.2 & 2.29 & 2.0 \\
        \fdgasoftg{0.9} & 2.93 & 12.6 & 2.71 & 2.0 \\
        \fdgasoftg{0.95} & 3.06 & 20.0 & 3.38 & 2.0 \\
        \fdgasoftg{0.99} & 3.35 & 60.6 & 6.49 & 3.0 \\
        \ldga & 3.31 & 100.0 & 12.35 & 3.0 \\
        \rdgasqrt & 3.21 & 100.0 & 12.11 & 4.0 \\\bottomrule
    \end{tabular}}{All values are averaged over five different instances.
    The values are computed for the set of orders in the \emph{last} phase of each run.}
  \end{table}

  \begin{figure}[H]
    \FIGURE{
    \includegraphics[width=0.8\textwidth]{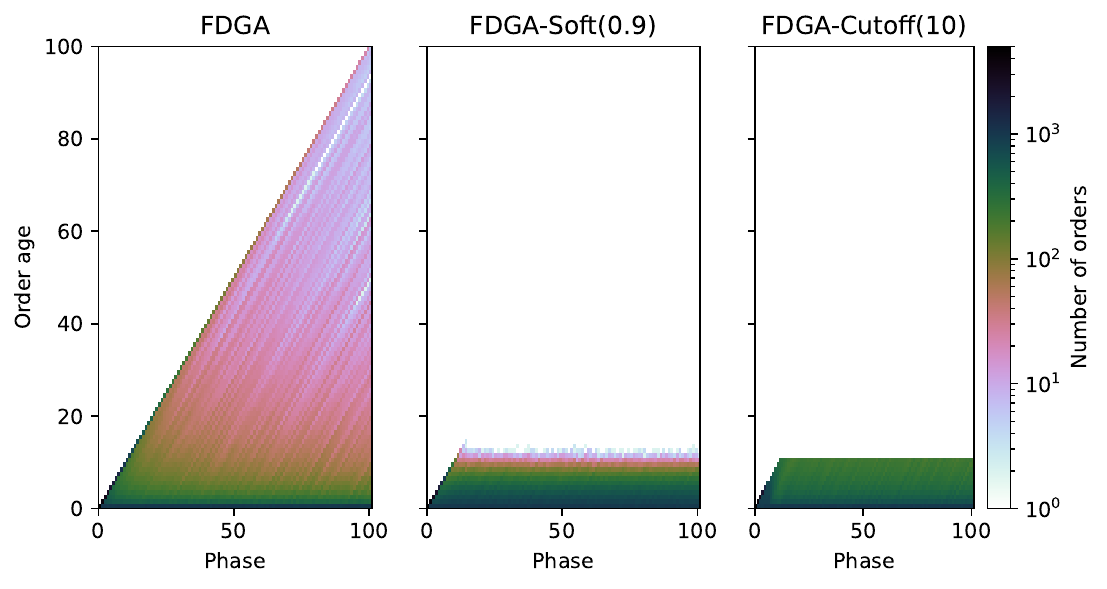}}
    {Age distribution data of single runs of \fdga, \fdgasoftg{0.9}, and \fdgacutoffn{10} on a medium-fewzones instance.
    \label{fig:app:fdga-age-fewzones}}{  While many orders are completed right away, many of them remain in the problem until the end of the experiment.
    }
  \end{figure}

  \begin{table}[H]
    \TABLE{Benchmark data comparing different multi-phase algorithms for medium-fewerzones instances.  \label{tab:app:multiphase-performance-medium-fewerzones}}
    {
      \begin{tabular}{c@{\hskip 2\tabcolsep}|cccc}\toprule
        Algorithm & Total cost & Rel.\ cost & \shortstack{Avg.\ cost\\ per item} & \shortstack{Avg.\ cost\\ per order} \\\midrule
        \fdga & 3745953 & (reference) & 14.36 & 37.38 \\
        \fdgacutoffn{5} & 4883920 & $+30.4\% $ & 18.98 & 49.99 \\
        \fdgacutoffn{10} & 4546505 & $+21.4\% $ & 17.63 & 46.2 \\
        \fdgacutoffn{15} & 4329714 & $+15.6\% $ & 16.74 & 43.79 \\
        \fdgacutoffn{20} & 4192427 & $+11.9\% $ & 16.19 & 42.3 \\
        \fdgacutoffn{25} & 4097160 & $+9.4\% $ & 15.81 & 41.27 \\
        \fdgacutoffn{30} & 4041441 & $+7.9\% $ & 15.57 & 40.64 \\
        \fdgalinearn{10} & 4621480 & $+23.4\% $ & 17.92 & 46.99 \\
        \fdgalinearn{15} & 4367900 & $+16.6\% $ & 16.9 & 44.22 \\
        \fdgalinearn{20} & 4233545 & $+13.0\% $ & 16.35 & 42.74 \\
        \fdgalinearn{25} & 4137448 & $+10.5\% $ & 15.94 & 41.61 \\
        \fdgalinearn{30} & 4073021 & $+8.7\% $ & 15.68 & 40.92 \\
        \fdgasoftg{0.8} & 3583825 & $-4.3\% $ & 13.88 & 36.49 \\
        \fdgasoftg{0.9} & 3476141 & $-7.2\% $ & 13.45 & 35.28 \\
        \fdgasoftg{0.95} & 3495369 & $-6.7\% $ & 13.49 & 35.32 \\
        \fdgasoftg{0.99} & 3659150 & $-2.3\% $ & 14.07 & 36.71 \\
        \ldga & 3863173 & $+3.1\% $ & 14.88 & 38.9 \\
        \rdgasqrt & 4449204 & $+18.8\% $ & 17.14 & 44.82 \\\bottomrule
    \end{tabular}}{The parameter for the number of new orders is set to 1000 if applicable and corresponds approximately to the number of new orders in each phase. The column for relative cost compares to the reference of \fdga. The values are averaged over 5 instances and rounded to appropriate precision.}
  \end{table}

  \begin{table}[H]
    \TABLE{Further computational measures for different algorithms on medium-fewerzones instances.
    \label{tab:multiphaseinsights-medium-fewerzones}}
    {
      \begin{tabular}{c@{\hskip 2\tabcolsep}|cccc}\toprule
        Algorithm & \shortstack{Average\\ order length} & \shortstack{Maximum\\ order age} & \shortstack{Average\\ order age} & \shortstack{Median\\ order age} \\\midrule
        \fdga & 3.63 & 100.0 & 15.52 & 4.0 \\
        \fdgacutoffn{5} & 2.74 & 5.0 & 2.33 & 2.0 \\
        \fdgacutoffn{10} & 3.07 & 10.0 & 3.54 & 3.0 \\
        \fdgacutoffn{15} & 3.21 & 15.0 & 4.52 & 3.0 \\
        \fdgacutoffn{20} & 3.3 & 20.0 & 5.48 & 3.0 \\
        \fdgacutoffn{25} & 3.36 & 25.0 & 6.32 & 3.6 \\
        \fdgacutoffn{30} & 3.4 & 30.0 & 7.11 & 3.8 \\
        \fdgalinearn{10} & 3.04 & 10.0 & 3.36 & 3.0 \\
        \fdgalinearn{15} & 3.21 & 15.0 & 4.49 & 3.0 \\
        \fdgalinearn{20} & 3.29 & 20.0 & 5.42 & 3.0 \\
        \fdgalinearn{25} & 3.36 & 25.0 & 6.3 & 3.2 \\
        \fdgalinearn{30} & 3.39 & 30.0 & 7.07 & 3.8 \\
        \fdgasoftg{0.8} & 2.88 & 9.4 & 2.43 & 2.0 \\
        \fdgasoftg{0.9} & 3.02 & 13.4 & 2.91 & 2.0 \\
        \fdgasoftg{0.95} & 3.14 & 21.4 & 3.69 & 2.4 \\
        \fdgasoftg{0.99} & 3.42 & 68.8 & 7.27 & 3.0 \\
        \ldga & 3.24 & 100.0 & 11.73 & 3.0 \\
        \rdgasqrt & 3.24 & 100.0 & 11.63 & 3.6 \\\bottomrule
    \end{tabular}}{All values are averaged over five different instances.
      The values are computed for the set of orders in the \emph{last} phase of each run.
    }
  \end{table}

  \begin{figure}[H]
    \FIGURE{
    \includegraphics[width=0.8\textwidth]{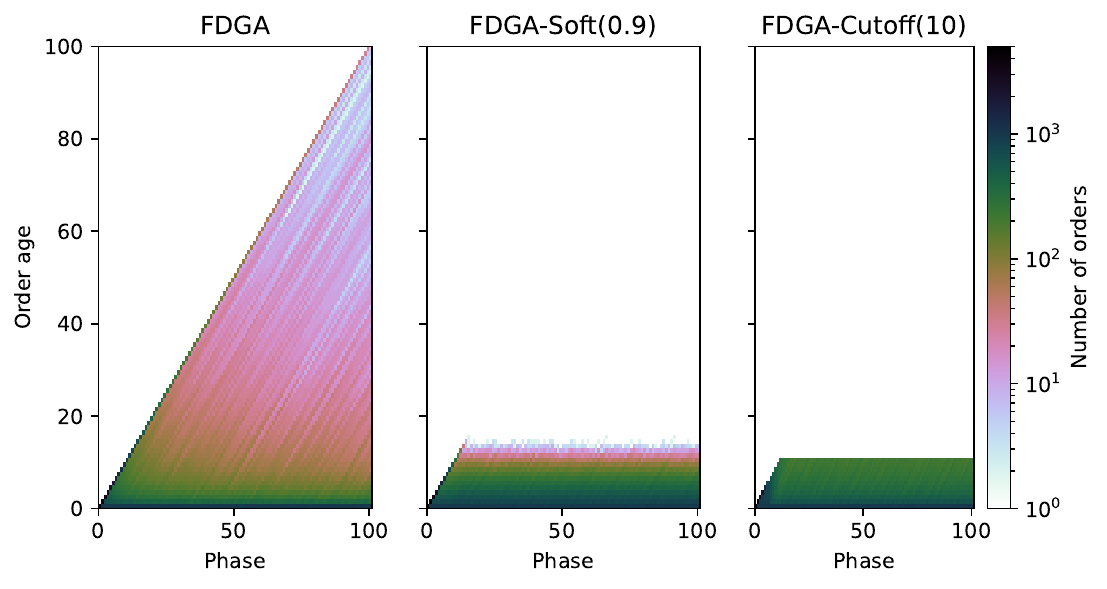}}
    {Age distribution data of single runs of \fdga, \fdgasoftg{0.9}, and \fdgacutoffn{10} on a medium-fewerzones instance.
    \label{fig:app:fdga-age-fewerzones}}{  While many orders are completed right away, many of them remain in the problem until the end of the experiment.
    }
  \end{figure}

  \begin{table}[H] 
    \TABLE{Benchmark data comparing different multi-phase algorithms for medium-longbatch instances.  \label{tab:app:multiphase-performance-medium-longbatch}}
    {
      \begin{tabular}{c@{\hskip 2\tabcolsep}|cccc}\toprule      Algorithm & Total cost & Rel.\ cost & \shortstack{Avg.\ cost\\ per item} & \shortstack{Avg.\ cost\\ per order} \\\midrule
        \fdga & 6563409 & (reference) & 24.61 & 64.41 \\
        \fdgacutoffn{5} & 8441069 & $+28.6\% $ & 31.62 & 83.37 \\
        \fdgacutoffn{10} & 7819809 & $+19.1\% $ & 29.17 & 76.66 \\
        \fdgacutoffn{15} & 7385350 & $+12.5\% $ & 27.62 & 72.51 \\
        \fdgacutoffn{20} & 7154438 & $+9.0\% $ & 26.77 & 70.21 \\
        \fdgacutoffn{25} & 7026946 & $+7.1\% $ & 26.3 & 68.96 \\
        \fdgacutoffn{30} & 6945899 & $+5.8\% $ & 26.0 & 68.16 \\
        \fdgalinearn{10} & 7871377 & $+19.9\% $ & 29.55 & 77.7 \\
        \fdgalinearn{15} & 7391335 & $+12.6\% $ & 27.76 & 72.86 \\
        \fdgalinearn{20} & 7147220 & $+8.9\% $ & 26.83 & 70.38 \\
        \fdgalinearn{25} & 6995720 & $+6.6\% $ & 26.29 & 68.92 \\
        \fdgalinearn{30} & 6910578 & $+5.3\% $ & 25.97 & 68.05 \\
        \fdgasoftg{0.8} & 6817485 & $+3.9\% $ & 25.46 & 67.13 \\
        \fdgasoftg{0.9} & 6500641 & $-1.0\% $ & 24.3 & 63.98 \\
        \fdgasoftg{0.95} & 6397924 & $-2.5\% $ & 23.95 & 62.97 \\
        \fdgasoftg{0.99} & 6468726 & $-1.4\% $ & 24.21 & 63.48 \\
        \ldga & 6412469 & $-2.3\% $ & 24.12 & 63.33 \\
        \rdgasqrt & 7052168 & $+7.4\% $ & 26.52 & 69.51 \\\bottomrule
      \end{tabular}
    }{The parameter for the number of new orders is set to 1000 if applicable and corresponds approximately to the number of new orders in each phase. The column for relative cost compares to the reference of \fdga. The values are averaged over 5 instances and rounded to appropriate precision.}
  \end{table}

  \begin{table}[H] 
    \TABLE{Further computational measures for different algorithms on medium-longbatch instances.
    \label{tab:multiphaseinsights-medium-longbatch}}
    {
      \begin{tabular}{c@{\hskip 2\tabcolsep}|cccc}\toprule
        Algorithm & \shortstack{Average\\ order length} & \shortstack{Maximum\\ order age} & \shortstack{Average\\ order age} & \shortstack{Median\\ order age} \\\midrule
        \fdga & 3.23 & 100.0 & 10.58 & 3.0 \\
        \fdgacutoffn{5} & 2.72 & 5.0 & 2.29 & 2.0 \\
        \fdgacutoffn{10} & 2.93 & 10.0 & 3.38 & 3.0 \\
        \fdgacutoffn{15} & 3.0 & 15.0 & 4.2 & 3.0 \\
        \fdgacutoffn{20} & 3.06 & 20.0 & 4.93 & 3.0 \\
        \fdgacutoffn{25} & 3.09 & 25.0 & 5.49 & 3.0 \\
        \fdgacutoffn{30} & 3.11 & 30.0 & 6.01 & 3.0 \\
        \fdgalinearn{10} & 2.89 & 10.0 & 3.25 & 3.0 \\
        \fdgalinearn{15} & 3.0 & 15.0 & 4.16 & 3.0 \\
        \fdgalinearn{20} & 3.05 & 20.0 & 4.87 & 3.0 \\
        \fdgalinearn{25} & 3.09 & 25.0 & 5.49 & 3.0 \\
        \fdgalinearn{30} & 3.12 & 30.0 & 6.04 & 3.0 \\
        \fdgasoftg{0.8} & 2.73 & 7.8 & 2.19 & 2.0 \\
        \fdgasoftg{0.9} & 2.81 & 11.2 & 2.52 & 2.0 \\
        \fdgasoftg{0.95} & 2.9 & 18.4 & 3.06 & 2.0 \\
        \fdgasoftg{0.99} & 3.08 & 53.8 & 5.14 & 3.0 \\
        \ldga & 2.97 & 100.0 & 8.27 & 3.0 \\
        \rdgasqrt & 3.11 & 100.0 & 10.15 & 3.0 \\\bottomrule
    \end{tabular}}{All values are averaged over five different instances.
      The values are computed for the set of orders in the \emph{last} phase of each run.
    }
  \end{table}

  \begin{figure}[H]
    \FIGURE{
    \includegraphics[width=0.8\textwidth]{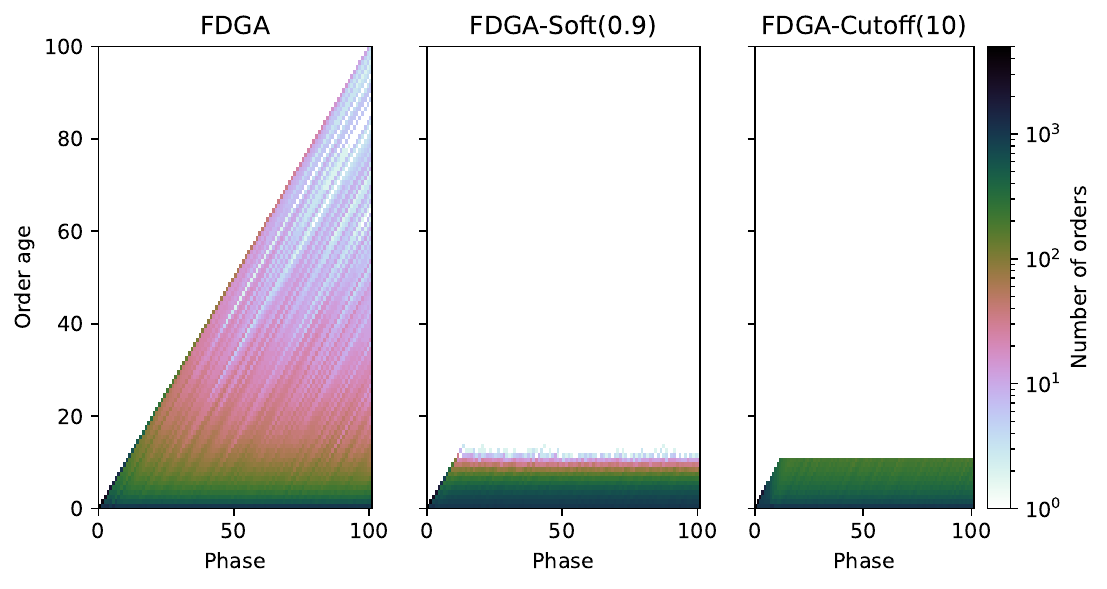}}
    {Age distribution data of single runs of \fdga, \fdgasoftg{0.9}, and \fdgacutoffn{10} on a medium-longbatch instance.
    \label{fig:app:fdga-age-longbatch}}{  While many orders are completed right away, many of them remain in the problem until the end of the experiment.
    }
  \end{figure}

  \begin{table}[H]
    \TABLE{Benchmark data comparing different multi-phase algorithms for large instances.  \label{tab:app:multiphase-performance-large}}
    {
      \begin{tabular}{c@{\hskip 2\tabcolsep}|cccc}\toprule      Algorithm & Total cost & Rel.\ cost & \shortstack{Avg.\ cost\\ per item} & \shortstack{Avg.\ cost\\ per order} \\\midrule
        \fdga & 87550699 & (reference) & 34.57 & 89.72 \\
        \ldga & 94800948 & $+8.3\% $ & 37.44 & 97.58 \\
        \rdgasqrt & 113283647 & $+29.4\% $ & 44.74 & 116.7 \\
        \fdgacutoffn{5} & 105116080 & $+20.1\% $ & 41.54 & 109.44 \\
        \fdgacutoffn{10} & 102316383 & $+16.9\% $ & 40.42 & 105.86 \\
        \fdgacutoffn{15} & 98299057 & $+12.3\% $ & 38.83 & 101.44 \\
        \fdgacutoffn{20} & 95859887 & $+9.5\% $ & 37.86 & 98.77 \\
        \fdgacutoffn{25} & 94138537 & $+7.5\% $ & 37.18 & 96.91 \\
        \fdgacutoffn{30} & 92982284 & $+6.2\% $ & 36.72 & 95.65 \\
        \fdgalinearn{10} & 103842705 & $+18.6\% $ & 41.03 & 107.54 \\
        \fdgalinearn{15} & 99300767 & $+13.4\% $ & 39.22 & 102.5 \\
        \fdgalinearn{20} & 96571183 & $+10.3\% $ & 38.14 & 99.51 \\
        \fdgalinearn{25} & 94728169 & $+8.2\% $ & 37.41 & 97.52 \\
        \fdgalinearn{30} & 93511512 & $+6.8\% $ & 36.93 & 96.2 \\
        \fdgasoftg{0.8} & 86994048 & $-0.6\% $ & 34.36 & 90.4 \\
        \fdgasoftg{0.9} & 83130094 & $-5.0\% $ & 32.83 & 86.18 \\
        \fdgasoftg{0.95} & 82055044 & $-6.3\% $ & 32.41 & 84.85 \\
        \fdgasoftg{0.99} & 84573403 & $-3.4\% $ & 33.39 & 86.99 \\\bottomrule
      \end{tabular}
    }{The parameter for the number of new orders is set to 10,000 if applicable and corresponds approximately to the number of new orders in each phase. The column for relative cost compares to the reference of \fdga. The values are averaged over 5 instances and rounded to appropriate precision.}
  \end{table}

  \begin{table}[H]
    \TABLE{Further computational measures for different algorithms on large instances.
    \label{tab:multiphaseinsights-large}}
    {
      \begin{tabular}{c@{\hskip 2\tabcolsep}|cccc}\toprule
        Algorithm & \shortstack{Average\\ order length} & \shortstack{Maximum\\ order age} & \shortstack{Average\\ order age} & \shortstack{Median\\ order age} \\\midrule
        \fdga & 3.86 & 100.0 & 23.32 & 8.0 \\
        \ldga & 3.45 & 100.0 & 17.37 & 5.0 \\
        \rdgasqrt & 3.37 & 100.0 & 18.12 & 6.0 \\
        \fdgacutoffn{5} & 2.72 & 5.0 & 2.36 & 2.0 \\
        \fdgacutoffn{10} & 3.07 & 10.0 & 3.73 & 3.0 \\
        \fdgacutoffn{15} & 3.26 & 15.0 & 4.96 & 4.0 \\
        \fdgacutoffn{20} & 3.37 & 20.0 & 6.11 & 4.0 \\
        \fdgacutoffn{25} & 3.45 & 25.0 & 7.21 & 4.0 \\
        \fdgacutoffn{30} & 3.51 & 30.0 & 8.25 & 5.0 \\
        \fdgalinearn{10} & 3.01 & 10.0 & 3.49 & 3.0 \\
        \fdgalinearn{15} & 3.24 & 15.0 & 4.82 & 4.0 \\
        \fdgalinearn{20} & 3.36 & 20.0 & 6.02 & 4.0 \\
        \fdgalinearn{25} & 3.44 & 25.0 & 7.12 & 4.0 \\
        \fdgalinearn{30} & 3.51 & 30.0 & 8.18 & 5.0 \\
        \fdgasoftg{0.8} & 2.79 & 8.2 & 2.29 & 2.0 \\
        \fdgasoftg{0.9} & 2.93 & 12.0 & 2.66 & 2.0 \\
        \fdgasoftg{0.95} & 3.1 & 19.0 & 3.37 & 2.4 \\
        \fdgasoftg{0.99} & 3.49 & 62.2 & 7.23 & 4.0 \\\bottomrule
    \end{tabular}}{All values are averaged over five different instances.
      The values are computed for the set of orders in the \emph{last} phase of each run.
    }
  \end{table}

  \begin{figure}[H]
    \FIGURE{
    \includegraphics[width=0.8\textwidth]{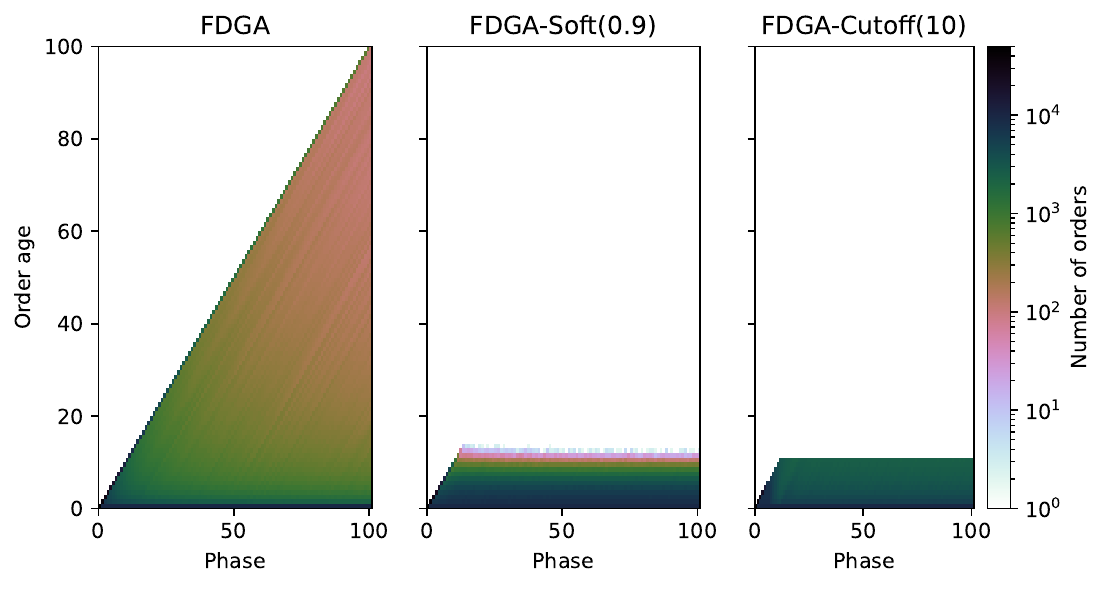}}
    {Age distribution data of single runs of \fdga, \fdgasoftg{0.9}, and \fdgacutoffn{10} on a large instance.
    \label{fig:app:fdga-age-large}}{While many orders are completed right away, many of them remain in the problem until the end of the experiment.
    }
  \end{figure}

  \begin{table}[H]
    \TABLE{Benchmark data comparing different multi-phase algorithms for large-fewzones instances.  \label{tab:app:multiphase-performance-large-fewzones}}
    {
      \begin{tabular}{c@{\hskip 2\tabcolsep}|cccc}\toprule      Algorithm & Total cost & Rel.\ cost & \shortstack{Avg.\ cost\\ per item} & \shortstack{Avg.\ cost\\ per order} \\\midrule
        \fdga & 43981250 & (reference) & 17.36 & 45.02 \\
        \fdgacutoffn{5} & 54296417 & $+23.5\% $ & 21.46 & 56.51 \\
        \fdgacutoffn{10} & 53246126 & $+21.1\% $ & 21.04 & 55.02 \\
        \fdgacutoffn{15} & 51400118 & $+16.9\% $ & 20.3 & 52.96 \\
        \fdgacutoffn{20} & 50092825 & $+13.9\% $ & 19.78 & 51.54 \\
        \fdgacutoffn{25} & 49028580 & $+11.5\% $ & 19.36 & 50.39 \\
        \fdgacutoffn{30} & 48185738 & $+9.6\% $ & 19.03 & 49.49 \\
        \fdgalinearn{10} & 54051806 & $+22.9\% $ & 21.36 & 55.9 \\
        \fdgalinearn{15} & 52250241 & $+18.8\% $ & 20.64 & 53.84 \\
        \fdgalinearn{20} & 50807093 & $+15.5\% $ & 20.07 & 52.28 \\
        \fdgalinearn{25} & 49638096 & $+12.9\% $ & 19.6 & 51.01 \\
        \fdgalinearn{30} & 48764224 & $+10.9\% $ & 19.26 & 50.09 \\
        \fdgasoftg{0.8} & 41342885 & $-6.0\% $ & 16.33 & 42.85 \\
        \fdgasoftg{0.9} & 40298227 & $-8.4\% $ & 15.92 & 41.66 \\
        \fdgasoftg{0.95} & 40632844 & $-7.6\% $ & 16.05 & 41.89 \\
        \fdgasoftg{0.99} & 42683060 & $-3.0\% $ & 16.85 & 43.81 \\
        \ldga & 51058898 & $+16.1\% $ & 20.16 & 52.61 \\
        \rdgasqrt & 57668601 & $+31.1\% $ & 22.77 & 59.32 \\\bottomrule
      \end{tabular}
    }{The parameter for the number of new orders is set to 10,000 if applicable and corresponds approximately to the number of new orders in each phase. The column for relative cost compares to the reference of \fdga. The values are averaged over 5 instances and rounded to appropriate precision.}
  \end{table}

  \begin{table}[H]
    \TABLE{Further computational measures for different algorithms on large-fewzones instances.
    \label{tab:multiphaseinsights-large-fewzones}}
    {
      \begin{tabular}{c@{\hskip 2\tabcolsep}|cccc}\toprule
        Algorithm & \shortstack{Average\\ order length} & \shortstack{Maximum\\ order age} & \shortstack{Average\\ order age} & \shortstack{Median\\ order age} \\\midrule
        \fdga & 3.97 & 100.0 & 21.71 & 7.0 \\
        \fdgacutoffn{5} & 2.73 & 5.0 & 2.36 & 2.0 \\
        \fdgacutoffn{10} & 3.16 & 10.0 & 3.68 & 3.0 \\
        \fdgacutoffn{15} & 3.36 & 15.0 & 4.88 & 3.6 \\
        \fdgacutoffn{20} & 3.48 & 20.0 & 6.02 & 4.0 \\
        \fdgacutoffn{25} & 3.56 & 25.0 & 7.09 & 4.0 \\
        \fdgacutoffn{30} & 3.62 & 30.0 & 8.08 & 4.8 \\
        \fdgalinearn{10} & 3.11 & 10.0 & 3.46 & 3.0 \\
        \fdgalinearn{15} & 3.35 & 15.0 & 4.77 & 3.0 \\
        \fdgalinearn{20} & 3.47 & 20.0 & 5.93 & 4.0 \\
        \fdgalinearn{25} & 3.56 & 25.0 & 7.02 & 4.0 \\
        \fdgalinearn{30} & 3.62 & 30.0 & 8.03 & 4.8 \\
        \fdgasoftg{0.8} & 2.93 & 9.8 & 2.48 & 2.0 \\
        \fdgasoftg{0.9} & 3.11 & 15.0 & 3.03 & 2.0 \\
        \fdgasoftg{0.95} & 3.3 & 23.6 & 3.94 & 3.0 \\
        \fdgasoftg{0.99} & 3.67 & 77.0 & 8.51 & 4.0 \\
        \ldga & 3.35 & 100.0 & 13.62 & 4.0 \\
        \rdgasqrt & 3.49 & 100.0 & 15.76 & 4.2 \\\bottomrule
    \end{tabular}}{All values are averaged over five different instances.
      The values are computed for the set of orders in the \emph{last} phase of each run.
    }
  \end{table}

  \begin{figure}[H]
    \FIGURE{
    \includegraphics[width=0.8\textwidth]{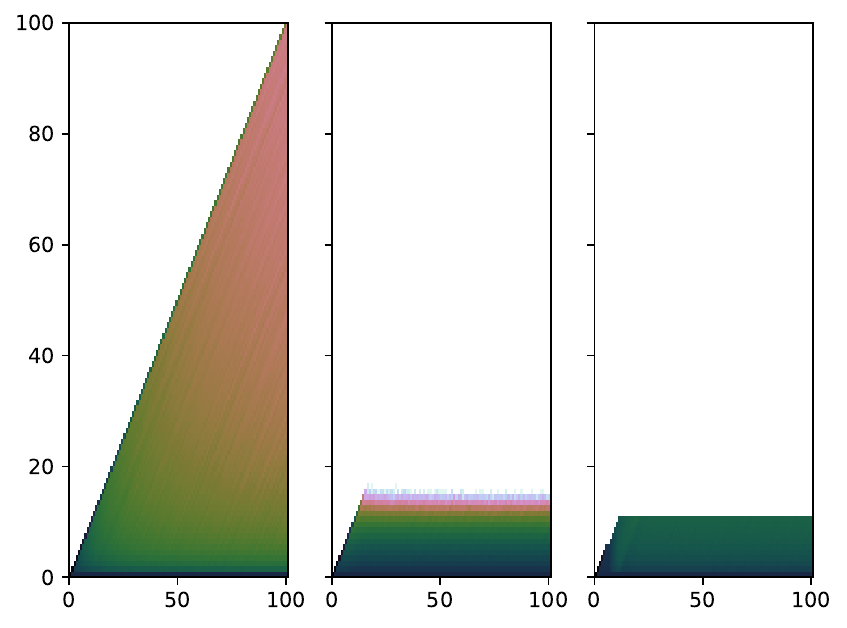}}
    {Age distribution data of single runs of \fdga, \fdgasoftg{0.9}, and \fdgacutoffn{10} on a large-fewzones instance.
    \label{fig:app:fdga-age-large-fewzones}}{While many orders are completed right away, many of them remain in the problem until the end of the experiment.
    }
  \end{figure}

\end{APPENDICES}

\end{document}